\documentclass[11pt,a4paper]{article}
\usepackage[T1]{fontenc}
\usepackage[utf8]{inputenc}
\IfFileExists{lmodern.sty}{\usepackage{lmodern}}{}
\usepackage{amsmath,amssymb,amsthm,mathtools,microtype}
\usepackage{authblk}
\usepackage[margin=26mm]{geometry}
\usepackage{enumitem,booktabs,xurl,longtable,array}
\usepackage{flafter}
\usepackage{needspace}
\usepackage{tikz}
\usepackage[normalem]{ulem}
\usetikzlibrary{positioning,fit,calc,arrows.meta,decorations.pathreplacing,decorations.pathmorphing,backgrounds,shapes.geometric}
\usepackage[hidelinks,unicode]{hyperref}
\usepackage{cleveref}
\usepackage{comment}

\newtheorem{theorem}{Theorem}[section]
\newtheorem{lemma}[theorem]{Lemma}

\newtheorem{corollary}[theorem]{Corollary}
\newtheorem{claim}{Claim}
\newtheorem{conjecture}{Conjecture}

\theoremstyle{definition}
\newtheorem{definition}[theorem]{Definition}
\theoremstyle{remark}
\newtheorem{remark}[theorem]{Remark}
\numberwithin{equation}{section}

\DeclareMathOperator{\ch}{ch}
\DeclareMathOperator{\st}{st}
\newcommand{\Pk}{\ensuremath{(\mathrm P_k)}}
\newcommand{\Mk}{\ensuremath{(\mathrm M_k)}}
\newcommand{\Van}{\mathrm{V}}

\newcommand{\RR}{\mathbb R}
\newcommand{\dcup}{\mathbin{\dot\cup}}

\setlist{itemsep=0.2em,topsep=0.4em}
\allowdisplaybreaks[1]
\tikzset{
  vtx/.style={circle,draw,fill=white,inner sep=1.5pt,minimum size=5.5pt},
  vtxS/.style={circle,draw,fill=black,inner sep=1.5pt,minimum size=5.5pt},
  clq/.style={draw,rounded corners=7pt,thick,inner sep=5pt},
  fib/.style={draw,rounded corners=4pt,inner sep=4pt,fill=gray!10,minimum height=9mm,minimum width=9mm,font=\small},
  cbox/.style={draw,rounded corners=5pt,inner sep=5pt,fill=gray!8,font=\small,align=center},
  nonedge/.style={densely dashed,thick,gray!60!black},
  full/.style={double,double distance=1.4pt,thick},
  incomp/.style={decorate,decoration={snake,amplitude=.5mm,segment length=2.5mm},thick},
  lab/.style={font=\small},
  stepnode/.style={draw,rounded corners=3pt,align=center,font=\small,inner sep=4pt,fill=gray!6},
  bad/.style={draw,rounded corners=3pt,align=center,font=\small,inner sep=4pt,fill=gray!25},
  arr/.style={-{Stealth[length=2.2mm]},thick}
}

\title{Claw-free perfect graphs are chromatic-choosable}
\author[1]{Suchismita Mishra}
\author[2]{Narayanan N}
\author[2]{Umesh Shankar}
\affil[1]{National Institute of Informatics, Tokyo}
\affil[2]{Department of Mathematics, Indian Institute of Technology Madras, Chennai, India}
\date{27 October 2026}
\hypersetup{pdftitle={Claw-free perfect graphs are chromatic-choosable},
pdfauthor={Suchismita Mishra; Narayanan N; Umesh Shankar},
pdfsubject={}}

\begin{document}
\maketitle
\begin{abstract}
We prove that every claw-free perfect graph $G$ satisfies $\ch(G)=\chi(G)=\omega(G)$.
More strongly, if $S$ is a clique of simplicial vertices of $G$ and $k\ge\omega(G)$,
every proper precoloring of $S$ extends to a proper coloring of $G$ from arbitrary
lists of size $k$ on $V(G)\setminus S$.  The proof has three coloring steps:
a mixed-list theorem for cobipartite graphs, a precoloring-extension theorem for
elementary graphs built on Galvin's kernel method and a homogeneous-pair
exchange argument, and a polynomial method (Combinatorial Nullstellensatz) argument for
graphs of interval type in which a Vandermonde factor is transferred across each
incomplete pair.  The three coloring steps are assembled along the structure theorems of
Chv\'atal--Sbihi, Maffray--Reed and Chudnovsky--Seymour, using double clique
cutsets for composition.  As a consequence, every perfect line graph is chromatic-choosable, recovering the theorem of Peterson and Woodall. Our precoloring extension theorem also resolves a conjecture posed by Chudnovsky and Seymour.

\end{abstract}
\textbf{Keywords :} List coloring, Claw-free perfect graphs,  List coloring conjecture.\\
\textbf{Mathematics Subject Classification :} 05C15, 05C17, 05C31.\\

\section{Introduction}

A central theme in list colouring is to determine when restrictions on the
available colours leave the chromatic number unchanged. A \emph{list
assignment} $L$ on a graph $G$ assigns to each vertex $v$ a set $L(v)$ of
admissible colours; an \emph{$L$-colouring} is a proper colouring $f$ with
$f(v)\in L(v)$ for all $v\in V(G)$. The \emph{choice number} $\ch(G)$ is the
least $k$ for which an $L$-colouring exists whenever $|L(v)|\ge k$ at every
vertex. Thus $\chi(G)\le\ch(G)$, and $G$ is \emph{chromatic-choosable} when
equality holds. List colouring was introduced by Vizing~\cite{Vizing} and
independently by Erd\H{o}s, Rubin and Taylor~\cite{ERT}. Even bipartite
graphs need not be chromatic-choosable, so perfectness alone does not force
$\ch(G)=\chi(G)$.

The List Colouring Conjecture asserts that every line graph is
chromatic-choosable. Galvin~\cite{Galvin} proved it for line graphs of
bipartite multigraphs, and Peterson and Woodall
\cite{PetersonWoodall,PetersonWoodallErratum} established it for perfect
line graphs. Since line graphs are claw-free, these results belong to the
broader conjecture of Gravier and Maffray~\cite{GravierMaffray} that
\emph{every claw-free graph is chromatic-choosable}. Recall that a graph is
claw-free if it contains no induced $K_{1,3}$, and perfect if each of its
induced subgraphs $H$ satisfies $\chi(H)=\omega(H)$.

The conjecture has been studied  intensively within claw-free
perfect graphs. Gravier and Maffray proved the result for
$3$-chromatic claw-free perfect graphs~\cite{GM2004}, extending their
earlier work on $3$-colourable elementary graphs~\cite{GM1997}.
Gravier, Maffray and Pastor~\cite{GMP} subsequently established it
whenever $\omega(G)\le 4$. More recently, Narayanan, Vasudevan and
Somasundaram~\cite{NNS} removed the clique-number restriction for the
entire class of elementary graphs. We settle the problem for all
claw-free perfect graphs.

\begin{theorem}\label{thm:main}
Every claw-free perfect graph $G$ satisfies
\[
  \ch(G)=\chi(G)=\omega(G).
\]
\end{theorem}

The difficulty in passing from basic structural classes to arbitrary claw-free
perfect graphs was identified by Chudnovsky and Seymour~\cite{CS}. Their
structure theorem uses compositions along cliques, but list colourings on
different pieces cannot in general be reconciled by permuting colour
names. Motivated by this obstacle, they formulated a precolouring-extension
conjecture for certain {\em simplicial cliques}.

\begin{conjecture}[Chudnovsky--Seymour {\cite[Conjecture~4.2]{CS}}]
Let $G$ be a perfect claw-free graph and let $K$ be a
simplicial clique in $G$. Let
$\{L_v\}_{v\in V(G)}$ be a set of lists of colours with
$|L_v|\geq \chi(G)$ for all $v\in V(G)$.
Let the vertices of $K$ be coloured in $|K|$ different
colours, each with a colour from its list.
Then this colouring of $K$ can be extended to a proper
list-colouring of $G$ with the lists
$\{L_v\}_{v\in V(G)}$.
\end{conjecture}

 We establish a stronger extension principle.
A vertex is \emph{simplicial} if its neighbourhood is a clique.

\begin{theorem}[Precolouring extension]\label{thm:main-P}
Let $G$ be a claw-free perfect graph, let $k\ge\omega(G)$ be an integer,
and let $S\subseteq V(G)$ be a clique of simplicial vertices, possibly
empty. If $\phi$ is any proper colouring of $G[S]$ and $L$ is any list
assignment on $V(G)\setminus S$ such that $|L(v)|\ge k$ for every
$v\notin S$, then $G$ has a proper colouring $f$ satisfying
\[
  f|_S=\phi
  \qquad\text{and}\qquad
  f(v)\in L(v)\quad(v\in V(G)\setminus S).
\]
\end{theorem}

Taking $S=\varnothing$ and $k=\omega(G)$ yields
Theorem~\ref{thm:main}.  Our precoloring extension theorem also resolves the above 
conjecture of Chudnovsky and Seymour. Their simplicial cliques are
cliques of simplicial vertices with the additional requirement that
$G-(S\cup N_G(S))$ be connected. Our theorem does not require this
connectivity condition. Moreover, the conjecture requires the prescribed
colours to belong to the original lists, whereas here they are arbitrary
and no lists on $S$ are needed.

The structural backbone for these results is the decomposition theorem
of Chv\'atal and Sbihi~\cite{ChvatalSbihi}. A \emph{clique cutset} of
a connected graph is a clique whose deletion disconnects the graph;
a connected graph without a clique cutset is an \emph{atom}.

\begin{theorem}[Chv\'atal--Sbihi~\cite{ChvatalSbihi}]\label{thm:atom}
Every claw-free perfect atom is elementary or peculiar.
\end{theorem}

Elementary graphs are those whose edges admit a two-colouring in which
the edges of every induced three-vertex path have different colours.
Maffray and Reed~\cite{MaffrayReed} characterized them as augmentations
of line graphs of bipartite multigraphs. Peculiar graphs have an
explicit clique-based description, and their chromatic-choosability
was established by Gravier and Maffray~\cite{GravierMaffray}.
Thus both kinds of atoms are known to be chromatic-choosable.
This fact does not by itself prove Theorem~\ref{thm:main}, however:
clique cutsets introduce compatibility constraints on lists, while
unrestricted clique sums need not preserve claw-freeness
\cite{ChudnovskyPlumettaz}. The extension theorem is designed to
meet precisely these difficulties.

\subsection{Proof strategy}\label{sec:strategy}

Our argument makes the structural decomposition compatible with lists by
alternating between precolouring extension and a \emph{mixed-list} property.
The mixed-list
property allows lists of size $k-r$ on an attachable clique $Q$ and
lists of size $k$ elsewhere, where $r\ge0$ and $|Q|+r\le k$.
If $S$ is a precoloured clique of simplicial vertices, then
$Q=N_G(S)$ is attachable in $G-S$; deleting the colours of $S$ from the
lists on $Q$ gives precisely such a mixed-list problem. This provides
the link between local colouring arguments and extension across
separators.

To prove the main result, three ingredients are needed. First, we prove an asymmetric mixed-list
theorem for cobipartite graphs. A saturation argument reduces the
problem to a canonical graph, where a nonzero coefficient of its graph
polynomial establishes the required choosability by the Combinatorial
Nullstellensatz of Alon~\cite{Alon}. Second, we establish an extension theorem
for elementary graphs. Galvin's kernel method handles the underlying
bipartite line graphs, while an extremal argument and a Hall-type
deficiency lemma control homogeneous pairs arising from incomplete
augmentations. These steps reduce a potential counterexample to the
line-graph case. Third, we treat graphs of interval type by organizing
their vertices into ordered clique fibres. The polynomial argument
eliminates these fibres successively, transferring and then recovering
Vandermonde factors across incomplete pairs. A refinement of the fibre
representation permits chosen simplicial vertices to be retained
throughout the elimination.

For the final step we use the structural refinement of Chudnovsky and
Seymour~\cite{CS}. Call a clique cutset \emph{rich} if each side of the
    separation contains two vertices that are not true twins. Their theorem
    states that a claw-free graph with a rich clique cutset is either of
    interval type or has a \emph{double clique cutset}: a clique cutset that
    splits into two cliques $C_1$, $C_2$ such that one side sees only $C_1$
    and the other only $C_2$ (precise definitions are in Section \ref{subsec:rich}). The
    extension property composes across double clique cutsets.

\section{Preliminaries}\label{sec:prelim}

\subsection{Notation and conventions}

We allow the empty graph and set its clique number, chromatic number
and choice number equal to zero. The empty set is a clique. Clique
parts in a cobipartition may be empty. All graph parameters and list
sizes used as bounds are integers. Empty products equal $1$, and
$0!=1$. In counterexamples to \Pk, we may assume without loss of
generality that every unprescribed list has size exactly $k$.

For $A\subseteq V(G)$, let $G[A]$ denote the subgraph induced by $A$,
and write $G-A=G[V(G)\setminus A]$. For a vertex $v$, let $N_G(v)$ and
$N_G[v]=N_G(v)\cup\{v\}$ denote its open and closed neighborhoods.
For a set $A$, its \emph{external neighborhood} is
\[
  N_G(A)=\left(\bigcup_{a\in A}N_G(a)\right)\setminus A.
\]
We omit the subscript when the graph is clear. Two adjacent vertices
are \emph{true twins} if their closed neighborhoods are equal.
Disjoint vertex sets form a \emph{complete pair} if every edge between
them is present, an \emph{anticomplete pair} if none is present, and an
\emph{incomplete pair} otherwise.

A graph is \emph{cobipartite} if its vertex set can be partitioned into
two cliques, and \emph{claw-free} if it has no induced subgraph
isomorphic to $K_{1,3}$. The \emph{line graph} of a graph $G$, denoted by
$\mathcal L(G)$, has vertex set $E(G)$, with two distinct vertices adjacent
if their corresponding edges have a common endpoint. We write $\st_G(w)$ for
the set of edges incident with a vertex $w$ of $G$.
  Multigraphs occur only as roots of line graphs; 
they may have parallel edges but no loops. 
For a list assignment $L$ and a set $T$ in its domain, put
$L(T)=\bigcup_{v\in T}L(v)$. For a coloring $f$, let
$f(T)=\{f(v):v\in T\}$. A proper coloring of a clique is injective,
so Hall's theorem gives
\begin{equation}\label{eq:hall}
  K\text{ is $L$-colorable}
  \quad\Longleftrightarrow\quad
  |L(T)|\ge |T|\quad\text{for every }T\subseteq K,
\end{equation}
whenever $K$ is a clique. The same inequalities characterize an
injective choice of colors on any vertex set, whether or not it is a
clique.

See Table~\ref{tab:notation} in the appendix for the recurring notation. Symbols introduced
inside an individual proof have the contextual meaning specified there.

\subsection{The two extension properties}

\begin{definition}\label{def:attachable}
A clique $K$ of a graph $G$ is \emph{attachable} if
$N_G(v)\setminus K$ is a clique for every $v\in K$. By convention, the
empty clique is attachable.
\end{definition}

We use repeatedly the following property which immediately follows from the definition of attachable clique.
\begin{equation}\label{eq:att-rule}
\text{If } x \in K \text{ has a neighbour $y \not\in K$ then every neighbour of $x$ that is not adjacent to $y$ lies in} K.
\end{equation}

\begin{definition}\label{def:PM}
Let $k\ge\omega(G)$ be an integer.
\begin{itemize}
\item A \emph{\Pk-instance} is a tuple $(G,S,\phi,L)$, where $S$ is a
clique of simplicial vertices of $G$, possibly empty, $\phi$ is a proper
coloring of $G[S]$, and $L$ is a list assignment on $V(G)\setminus S$
with $|L(v)|\ge k$ for every $v\notin S$. The instance is
\emph{colorable} if $G$ has a proper coloring extending $\phi$ and using
a color from $L(v)$ at every $v\notin S$. We say that $G$ satisfies
\Pk\ if every \Pk-instance on $G$ is colorable.

\item An \emph{\Mk-instance} is a tuple $(G,K,r,L)$, where $r\ge0$ is
an integer, $K$ is an attachable clique of $G$ with $|K|+r\le k$, and
$L$ is a list assignment on $V(G)$ satisfying
\[
  |L(v)|\ge
  \begin{cases}
    k-r,&v\in K,\\
    k,&v\notin K.
  \end{cases}
\]
The instance is \emph{colorable} if $G$ has a proper $L$-coloring.
We say that $G$ satisfies \Mk\ if every \Mk-instance on $G$ is colorable.
\end{itemize}
\end{definition}

The next lemma explains both the definition of attachability and the
form of the reduced list bound in \Mk.

\begin{lemma}[Simplicial cliques]\label{lem:simplicial-clique}
Let $K$ be a clique of simplicial vertices in a graph $G$, and put
$Q=N_G(K)$. Then the vertices of $K$ are pairwise true twins and
$K\cup Q$ is a clique. Moreover, if $G$ is claw-free, then:
\begin{enumerate}[label=\textup{(\alph*)}]
\item The clique $Q$ is attachable in $G-K$.
\item Given a \Pk-instance $(G,K,\phi,L)$, define lists on $G-K$ by
\[
  L'(v):=
  \begin{cases}
    L(v)\setminus\phi(K),&v\in Q,\\
    L(v),&v\notin K\cup Q.
  \end{cases}
\]
Then $(G-K,Q,|K|,L')$ is an \Mk-instance, and it is colorable if and
only if $(G,K,\phi,L)$ is colorable.
\end{enumerate}
\end{lemma}

\begin{proof}
The assertions are immediate when $K=\varnothing$, so assume that
$K\ne\varnothing$. If $u,v\in K$ are distinct, simpliciality of $u$
implies that every neighbor of $u$ other than $v$ is adjacent to $v$.
Thus $N_G[u]\subseteq N_G[v]$, and the reverse inclusion follows
symmetrically. Consequently all vertices of $K$ have closed
neighborhood $K\cup Q$, which is a clique.

Suppose that $G$ is claw-free. If a vertex $q\in Q$ had nonadjacent
neighbors $x,y$ outside $K\cup Q$, then, for any $u\in K$, the vertices
$q,u,x,y$ would induce a claw centered at $q$. Thus $Q$ is attachable
in $G-K$, proving (a).

For (b), the clique $K\cup Q$ gives $|K|+|Q|\le\omega(G)\le k$.
At most $|K|$ colors are removed from each list on $Q$, so
\[
  |L'(v)|\ge
  \begin{cases}
    k-|K|,&v\in Q,\\
    k,&v\in V(G)\setminus(K\cup Q).
  \end{cases}
\]
Together with (a), this makes $(G-K,Q,|K|,L')$ an \Mk-instance.
Every proper coloring of $G$ extending $\phi$ and respecting $L$
outside $K$ restricts to an $L'$-coloring of $G-K$, since $Q$ is
complete to $K$. Conversely, an $L'$-coloring
of $G-K$ combines with $\phi$ to color $G$, because $Q=N_G(K)$ and
no vertex of $Q$ uses a color in $\phi(K)$.
\end{proof}


\begin{lemma}[Minimum counterexamples]\label{lem:minimal}
Fix a nonnegative integer $k$ and a class $\mathcal C$ closed under
induced subgraphs. Among all uncolorable \Pk-instances
$(G,S,\phi,L)$ with $G\in\mathcal C$ and $\omega(G)\le k$, choose one
with $|V(G)|$ minimum. Then $G$ is connected, every simplicial vertex
of $G$ belongs to $S$, and every proper induced subgraph of $G$
satisfies \Pk.
\end{lemma}

\begin{proof}
Minimality is taken over all permissible choices of $S,\phi,L$.
Consequently, every proper induced subgraph satisfies \Pk\ for every
permissible instance on that subgraph. In particular, this is not a
statement restricted to the inherited lists.

The empty graph is colorable. If $G$ were disconnected, each component
would be a proper induced subgraph. A nonempty clique $S$ lies in one
component; in that component use its prescribed coloring, and in every
other component use the empty prescribed clique. The resulting
colorings combine to color $G$, a contradiction.

Suppose that a simplicial vertex $v$ lies outside $S$. In $G-v$ the
vertices of $S$ remain simplicial, so minimality gives a coloring
extending $\phi$. The set $N_G[v]$ is a clique, whence
$|N_G(v)|\le\omega(G)-1\le k-1$. At least one of the $k$ or more
colors in $L(v)$ is absent from its neighborhood. Giving that color to
$v$ completes a coloring of $G$, again a contradiction.
\end{proof}

\subsection{Minimum order and constrained edge maximality}\label{sec:extremal-setup}

Every minimum-counterexample argument fixes the integer $k$ first.
Minimum order is taken over \emph{all} uncolorable \Pk-instances in
the specified induced-subgraph-closed class $\mathcal C$ with clique
number at most $k$. Thus it is available for every permissible list
assignment on a smaller graph, including one obtained by changing lists.
Lemma~\ref{lem:minimal} records its immediate consequences.

Some arguments add a second extremal condition: with the vertex set, $S$, $\phi$ and the lists fixed (lists trimmed to size exactly $k$ as in Section 2.1), the number of edges is maximized among counterexamples in $\mathcal C$ with clique number at most $k$ in which $S$ is still a clique of simplicial vertices. Such a choice exists because there are finitely many graphs on a fixed vertex set, adding edges cannot create a coloring, and the instance still has minimum order. Whenever an edge is added in a proof, membership in $\mathcal C$, the clique bound and simpliciality of $S$ are checked before edge maximality is invoked.

Theorem~\ref{thm:M} uses a different, explicitly stated saturation:
the vertex set, lists, and designated attachable clique are fixed,
and only cobipartite supergraphs preserving attachability and the clique
bound are considered. No minimum-order assumption is used there.
The interval-type and final assembly arguments need minimum order only.

All extremal choices are local to the proof in which they are made.
An extremal property is used only under the hypotheses specified
for that choice. In particular, when an earlier theorem is applied
later, only its stated conclusion is used; the counterexample
selected in the later proof need not satisfy the extremal
conditions employed in proving the earlier theorem.

\subsection{Cobipartite graphs and matching saturation}

Let $H$ be a cobipartite graph with clique partition
$V(H)=X\dcup Y$. Its \emph{nonedge graph} $F$ is the graph on the
same vertex set whose edges are the nonadjacent $X$--$Y$ pairs of
$H$. Thus $F=\overline H$ is bipartite. A set is a clique of $H$ if
and only if it is independent in $F$. By K\H{o}nig's theorem, the
minimum size of a vertex cover of $F$ equals its matching number
$\nu(F)$; therefore
\begin{equation}\label{eq:konig}
  \omega(H)=|X|+|Y|-\nu(F).
\end{equation}

\begin{lemma}[Saturation]\label{lem:saturation}
Let $M$ be a maximum matching of $F$, and put
$H^+=K_{X\cup Y}-M$. Then $H\subseteq H^+$ and
$\omega(H^+)=\omega(H)$. In particular, every $L$-coloring of $H^+$
is an $L$-coloring of $H$.
\end{lemma}

\begin{proof}
Only nonedges of $H$ are added, and the nonedge graph of $H^+$ is the
matching $M$. Its matching number is $|M|=\nu(F)$, so
\eqref{eq:konig} gives the equality of clique numbers. A coloring
proper on a supergraph is proper on the original graph.
\end{proof}

Pairing the ends of $M$ and giving every remaining vertex its own
color also shows that $\chi(H)=\omega(H)$. The argument applies to
each induced subgraph, so cobipartite graphs are perfect. They are
claw-free because they have no independent set of size three.

Saturation preserves the clique number, but it need not preserve a
designated attachable clique. The proof of Theorem~\ref{thm:M} will
therefore add edges only when attachability is preserved as well.

\subsection{A coefficient criterion for list coloring}

Fix a linear order on $V(G)$. The \emph{graph polynomial} is
\[
  P_G(\mathbf{x})=
  \prod_{\substack{uv\in E(G)\\u<v}}(x_v-x_u)
  \in\mathbb R[x_v:v\in V(G)].
\]
For nonnegative integers $d(v)$, write
$\mathbf{x}^{d}=\prod_v x_v^{d(v)}$ and let
$[\mathbf{x}^{d}]P$ denote its coefficient in $P$.
The following interpolation identity yields the form of the
Combinatorial Nullstellensatz~\cite{Alon} used below.

\begin{lemma}[Coefficient formula]\label{lem:coefficient-formula}
Let $P\in\mathbb R[x_1,\ldots,x_n]$, and let
$d_1,\ldots,d_n$ be nonnegative integers with
$\deg P\le\sum_i d_i$. For each $i$, let $T_i\subseteq\mathbb R$
have $d_i+1$ elements and put
$\psi_i(x)=\prod_{a\in T_i}(x-a)$. Then
\begin{equation}\label{eq:coefficient-formula}
  [x_1^{d_1}\cdots x_n^{d_n}]P
  =\sum_{(a_1,\ldots,a_n)\in T_1\times\cdots\times T_n}
    \frac{P(a_1,\ldots,a_n)}{\prod_{i=1}^n\psi_i'(a_i)}.
\end{equation}
In particular, a nonzero coefficient implies that $P$ is nonzero at
some point of this Cartesian product.
\end{lemma}

\begin{proof}
For a set $T$ of $d+1$ distinct real numbers, Lagrange interpolation
gives
\[
  \sum_{a\in T}\frac{a^j}{\prod_{b\in T\setminus\{a\}}(a-b)}
  =\begin{cases}0,&0\le j<d,\\1,&j=d.\end{cases}
\]
All denominators are nonzero. Apply the product of these linear
functionals to each monomial $x_1^{e_1}\cdots x_n^{e_n}$ of $P$.
If $e_i<d_i$ for some $i$, the result is zero. If $e_i>d_i$ for some
$i$, the total-degree bound forces $e_j<d_j$ for another $j$, so the
result is again zero. The only surviving monomial has $e_i=d_i$ for
all $i$, and its multiplier is $1$. Linearity proves the identity.
For $n=0$, it is the identity for a constant polynomial, with the
usual empty-product convention.
\end{proof}

All colors appearing in a finite list-coloring instance, including
prescribed colors, can be identified injectively with real numbers.
We make this identification whenever graph polynomials are used.

\begin{corollary}\label{cor:nullstellensatz}
Let $Z\subseteq V(G)$, and let $\phi$ be a \emph{proper coloring of
$G[Z]$}. Define the polynomial in the variables indexed by
$V(G)\setminus Z$ by
\[
  P_{G,\phi}=
  \prod_{\substack{uv\in E(G-Z)\\u<v}}(x_v-x_u)
  \prod_{\substack{vz\in E(G)\\v\notin Z,\ z\in Z}}
       (x_v-\phi(z)).
\]
Suppose that nonnegative integers $d(v)$, $v\notin Z$, satisfy
\[
  \sum_{v\notin Z}d(v)=\deg P_{G,\phi}
  \qquad\text{and}\qquad
  \left[\prod_{v\notin Z}x_v^{d(v)}\right]P_{G,\phi}\ne0.
\]
Then $\phi$ extends to a proper coloring of $G$ from every list
assignment on $V(G)\setminus Z$ with $|L(v)|\ge d(v)+1$.
\end{corollary}

\begin{proof}
For each $v\notin Z$, choose a subset of $L(v)$ of size $d(v)+1$ and
use it as the interpolation set in
Lemma~\ref{lem:coefficient-formula}. At a point where $P_{G,\phi}$
is nonzero, adjacent uncolored vertices receive different colors, and
each uncolored vertex avoids the colors on its neighbors in $Z$.
Edges within $Z$ are already properly colored by hypothesis. The
resulting coloring therefore extends $\phi$ and respects every list.
If $Z=V(G)$, the polynomial is $1$ and the conclusion simply says
that the already proper coloring requires no extension.
\end{proof}

\section{Cobipartite graphs: the mixed-list theorem}\label{sec:cobipartite}

We prove that in a cobipartite graph an attachable clique $C$ can be
colored from lists of size $|C|$, provided that all remaining lists
have size at least the clique number. The short lists need not have
any relation to the other lists. The proof has two parts: an explicit
coefficient calculation for a canonical graph, and a saturation
argument reducing a hypothetical counterexample to that graph.

\subsection{Asymmetric lists}

We first derive an asymmetric list-coloring lemma from~\cite[Corollary~2.6(b)]{NNS}. Besides handling several boundary cases, it
will be used again in the elementary-graph argument.

\begin{lemma}[Asymmetric cobipartite lists]\label{lem:cobi-asymmetric}
Let $H$ admit a clique partition $V(H)=X\dcup Y$, where $|X|=a$ and $|Y|=b$. If a list assignment $L$ satisfies 
\[
  |L(x)|\ge \omega(H)\ \text{ for all } x\in X \qquad\text{and}\qquad |L(y)|\ge b\ \text{ for all } y\in Y
\]
then $H$ is $L$-colorable.
\end{lemma}

\begin{proof}
The assertion is immediate if either part is empty.
Otherwise, let $M$ be a maximum matching in the bipartite
graph of missing $X$--$Y$ pairs. By
Lemma~\ref{lem:saturation}, adding all missing pairs outside
$M$ yields a supergraph
\[
H^+=K_{X\cup Y}-M
\]
with $\omega(H^+)=\omega(H)=q$.
The lists have size at least $q\ge a$ on $X$ and at least
$b$ on $Y$. Corollary~2.6 (b) of~\cite{NNS} therefore gives an
$L$-coloring of $H^+$, which also colors $H$.
\end{proof}

\subsection{The canonical graph and its coefficient}

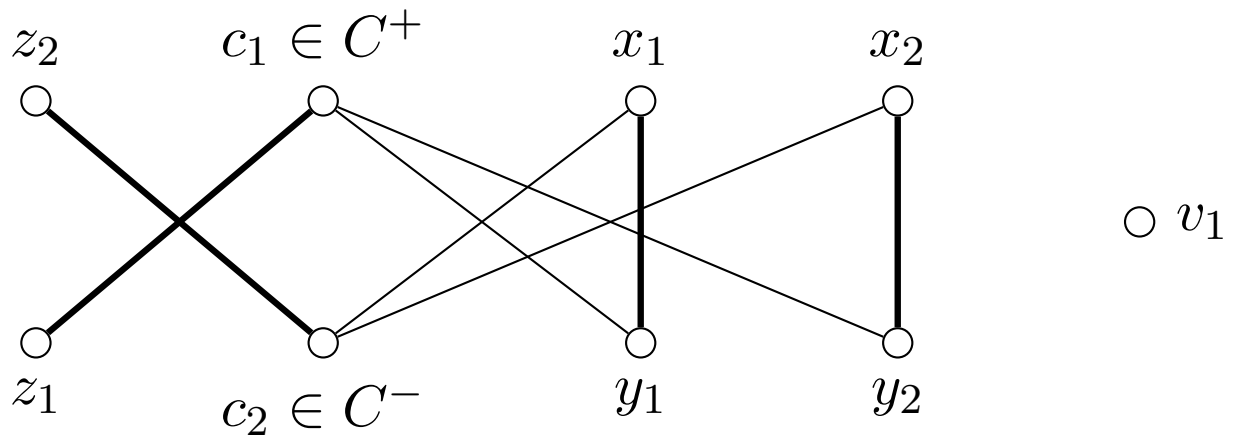
\begin{figure}[!htbp]
\centering
\begin{tikzpicture}
\node[vtx,label=above:$z_2$] (z2) at (0,1.6) {};
\node[vtx,label=above:{$c_1\in C^+$}] (c1) at (1.9,1.6) {};
\node[vtx,label=above:$x_1$] (x1) at (4.0,1.6) {};
\node[vtx,label=above:$x_2$] (x2) at (5.7,1.6) {};
\node[vtx,label=below:$z_1$] (z1) at (0,0) {};
\node[vtx,label=below:{$c_2\in C^-$}] (c2) at (1.9,0) {};
\node[vtx,label=below:$y_1$] (y1) at (4.0,0) {};
\node[vtx,label=below:$y_2$] (y2) at (5.7,0) {};
\node[vtx,label=right:$v_1$] (v1) at (7.3,.8) {};
\draw[very thick] (c1)--(z1);
\draw[very thick] (c2)--(z2);
\draw[very thick] (x1)--(y1);
\draw[very thick] (x2)--(y2);
\draw (c1)--(y1); \draw (c1)--(y2);
\draw (c2)--(x1); \draw (c2)--(x2);
\end{tikzpicture}
\caption{The nonedge graph of
$J(2,2,1;\{c_1\},\{c_2\})$. Thick edges are the matching pairs
$c_i z_i$ and $x_j y_j$; thin edges are the additional pairs in
\eqref{eq:J-nonedges}. The isolated vertex $v_1$ is universal in
$J$. Two distinct vertices are adjacent in $J$ exactly when they
are not joined in this drawing.}
\label{fig:canonical}
\end{figure}

Let $m,t\ge1$ and $u\ge0$ be integers, and put
$q_0=m+t$ and $q=q_0+u$. Start with pairwise disjoint sets
\[
\begin{gathered}
 C=\{c_1,\ldots,c_m\},\qquad Z=\{z_1,\ldots,z_m\},\qquad
 U=\{v_1,\ldots,v_u\},\\
 X=\{x_1,\ldots,x_t\},\qquad Y=\{y_1,\ldots,y_t\},
\end{gathered}
\]
and a partition $C=C^+\dcup C^-$.  Either of the two parts may be empty. Now, define
$J=J(m,t,u;C^+,C^-)$ by specifying that its nonedges are exactly
\begin{equation}\label{eq:J-nonedges}
\begin{gathered}
 c_i z_i\quad(1\le i\le m),\qquad
 x_j y_j\quad(1\le j\le t),\\
 c y_j\quad(c\in C^+,\ 1\le j\le t),\qquad
 c x_j\quad(c\in C^-,\ 1\le j\le t).
\end{gathered}
\end{equation}
This means that every vertex of $U$ is universal in $J$. 
Figure~\ref{fig:canonical} illustrates the nonedge graph.

Now, the  complement of $J-U$ has bipartition
\[
  C^+\cup X\cup\{z_i:c_i\in C^-\},
  \qquad
  C^-\cup Y\cup\{z_i:c_i\in C^+\}.
\]
Both parts have size $q_0$. The vertices of $U$ are isolated in the
nonedge graph and can be assigned to either part, so $J$ is
cobipartite. Its nonedge graph has a matching of size $q_0$, consisting
of the pairs $c_i z_i$ and $x_j y_j$. No larger matching is possible,
since there are only $2q_0$ nonisolated vertices. Thus
\eqref{eq:konig} gives
\[
  \omega(J)=2q_0+u-q_0=q.
\]
The clique $C$ is attachable: outside $C$, the neighbors of $c_i$
are $Z\setminus\{z_i\}$, all of $U$, and all of $X$ or all of $Y$
according as $c_i\in C^+$ or $c_i\in C^-$. These sets together form
a clique.

The exponents below deliberately give every vertex of $C$ the
smallest bound compatible with coloring the clique $C$. The other
exponents are chosen so that an interpolation grid forces colors to
occur on matched nonedges. This makes the coefficient calculation
explicit.

Order the vertices as
\[
 c_1,z_1,\ldots,c_m,z_m,x_1,y_1,\ldots,x_t,y_t,v_1,\ldots,v_u,
\]
and define $P_J$ using this order.

\begin{lemma}[Canonical coefficient]\label{lem:canonical-coefficient}
Define nonnegative integers $d(w)$ by
\[
\begin{aligned}
 d(c_i)&=m-1, & d(z_i)&=q-1,\\
 d(x_j)&=d(y_j)=q-1, & d(v_\ell)&=m+\ell-1.
\end{aligned}
\]
Then $\sum_w d(w)=|E(J)|$ and
\begin{equation}\label{eq:canonical-coefficient}
  [\mathbf{x}^{d}]P_J=(-1)^{m(t+u)}m!\,t!\ne0.
\end{equation}
Consequently, $J$ is colorable from arbitrary lists having at least
$m$ colors on $C$, at least $q$ on $Z\cup X\cup Y$, and at least
$m+\ell$ at $v_\ell$ for $1\le\ell\le u$.
\end{lemma}

\begin{proof}
We first verify the degree condition, then identify every nonzero
term in the coefficient formula and compute its value.

\emph{Degree.}
Put $J_0=J-U$. This graph has $2q_0$ vertices and $q_0+mt$
nonedges. Hence
\[
\begin{aligned}
  |E(J_0)|
    &=\binom{2q_0}{2}-q_0-mt\\
    &=m(m-1)+(m+2t)(q_0-1),
\end{aligned}
\]
which is the exponent sum when $u=0$. Adjoining $U$ adds
$2q_0u+\binom u2$ edges. The exponent sum increases by
\[
  u(m+2t)+\sum_{\ell=1}^u(m+\ell-1)
       =2q_0u+\binom u2.
\]
Thus the exponent sum equals $|E(J)|$ for every $u\ge0$.

\emph{Interpolation sets.}
In Lemma~\ref{lem:coefficient-formula}, use the sets
\[
\begin{aligned}
 T(c_i)&=\{1,\ldots,m\},\\
 T(z_i)=T(x_j)=T(y_j)&=\{1,\ldots,q\},\\
 T(v_\ell)&=\{1,\ldots,m\}
                  \cup\{q_0+1,\ldots,q_0+\ell\}.
\end{aligned}
\]
Each has $d(w)+1$ elements. Let $a$ be a point of this grid with
$P_J(a)\ne0$, viewed as a proper coloring of $J$. The clique $C$
uses all colors $1,\ldots,m$. Since $v_1$ is universal, its only
available color is $h_1=q_0+1$. Inductively, $v_\ell$ must use
$h_\ell=q_0+\ell$: all its other available colors occur on $C$ or
on $v_1,\ldots,v_{\ell-1}$.

Every vertex of $J_0$ is adjacent to all of $U$, so $J_0$ uses only
the colors $1,\ldots,q_0$. Its two clique parts both have $q_0$
vertices. Each color therefore occurs once in each part, and its two
occurrences form a nonedge of $J_0$. The only nonneighbor of $z_i$
in $J_0$ is $c_i$, so $a(z_i)=a(c_i)$. These pairs use the colors
$1,\ldots,m$. The remaining colors $m+1,\ldots,q_0$ must occur on
the pairs $x_j,y_j$. Conversely, all such assignments are proper.
Thus there are exactly $m!t!$ nonzero grid points: independently
permute the first $m$ colors among the $c_i,z_i$ pairs and the last
$t$ colors among the $x_j,y_j$ pairs.

\emph{The value of a summand for $u=0$.}
Write $A_0=\{1,\ldots,m\}$ and $B_0=\{m+1,\ldots,q_0\}$, and put
\[
\begin{gathered}
 V_A=\prod_{\substack{a<a'\\a,a'\in A_0}}(a'-a),\qquad
 V_B=\prod_{\substack{b<b'\\b,b'\in B_0}}(b'-b),\\
 \Gamma=\prod_{a\in A_0}\prod_{b\in B_0}(b-a).
\end{gathered}
\]
Start with the complete graph on $V(J_0)$ minus the $q_0$
equal-colored pairs. Two such pairs of distinct colors contribute
the fourth power of their color difference: their vertices form
consecutive two-vertex blocks in our order. To obtain $J_0$, remove
one further edge between each $C$-pair and each $X,Y$-pair. Its
factor is always $b-a$, where $a\in A_0$ and $b\in B_0$, since the
$C$-pair precedes the $X,Y$-pair. It follows that
\begin{equation}\label{eq:grid-numerator}
  P_{J_0}(a)=V_A^4 V_B^4\Gamma^3.
\end{equation}
In particular, this numerator is independent of both permutations
and of the partition $C^+,C^-$.

Let $\varphi_s(x)=\prod_{j=1}^s(x-j)$. The denominator in the
coefficient formula is
\[
  D_0=\prod_{a\in A_0}\varphi_m'(a)\varphi_{q_0}'(a)
       \prod_{b\in B_0}\varphi_{q_0}'(b)^2.
\]
For $a\in A_0$ and $b\in B_0$,
\[
\begin{aligned}
 \varphi_{q_0}'(a)
   &=\varphi_m'(a)\prod_{b\in B_0}(a-b),\\
 \varphi_{q_0}'(b)
   &=\prod_{a\in A_0}(b-a)
       \prod_{b'\in B_0\setminus\{b\}}(b-b').
\end{aligned}
\]
The products of the squared within-set derivative factors are
$V_A^4$ and $V_B^4$. The cross factors contribute
$(-1)^{mt}\Gamma$ from the first line and $\Gamma^2$ from the
second. Consequently,
\begin{equation}\label{eq:grid-denominator}
  D_0=(-1)^{mt}V_A^4V_B^4\Gamma^3.
\end{equation}
Every nonzero summand therefore equals $(-1)^{mt}$.

\emph{Adjoining universal vertices.}
Suppose $u\ge1$, and compare the calculation with the one for $u-1$.
Write $h=h_u=q_0+u$. Every interpolation set on $Z\cup X\cup Y$
gains the element $h$; the sets on $C$ and on the earlier universal
vertices are unchanged. The numerator gains the factor
\[
  \prod_{a\in A_0}(h-a)^2
  \prod_{b\in B_0}(h-b)^2
  \prod_{\ell<u}(h-h_\ell).
\]
The denominator gains
\[
  \left(\prod_{a\in A_0}(a-h)
         \prod_{b\in B_0}(b-h)^2\right)
  \left(\prod_{a\in A_0}(h-a)
         \prod_{\ell<u}(h-h_\ell)\right).
\]
The second parenthesis is the derivative factor for the new vertex
$v_u$. Taking the quotient shows that every summand is multiplied
by
\[
  \prod_{a\in A_0}\frac{h-a}{a-h}=(-1)^m.
\]
After all $u$ universal vertices have been adjoined, each of the
$m!t!$ nonzero summands is $(-1)^{m(t+u)}$. Summing proves
\eqref{eq:canonical-coefficient}. The list-coloring conclusion
follows from Corollary~\ref{cor:nullstellensatz} with
$Z=\varnothing$.
\end{proof}

\subsection{The mixed-list theorem}

\begin{theorem}[Mixed lists on cobipartite graphs]\label{thm:M}
Let $H$ be a cobipartite graph, let $C$ be an attachable clique of
$H$, and let $k\ge\omega(H)$ be an integer. Then $H$ is colorable
from every list assignment satisfying
\[
  |L(c)|\ge |C|\quad(c\in C),
  \qquad
  |L(v)|\ge k\quad(v\notin C).
\]
\end{theorem}

\begin{proof}
\setcounter{claim}{0}
If $C=\varnothing$, take any clique partition $X,Y$ of $H$ and
apply Lemma~\ref{lem:cobi-asymmetric}; the lists have size at least
$k\ge\omega(H)\ge |Y|$. This includes the empty graph. If
$C=V(H)$, Hall's theorem colors the clique $C$ from its lists.

Suppose, for a contradiction, that the theorem fails for an instance
$(H,C,L,k)$ not covered by these cases. Put $A=V(H)\setminus C$ and
$m=|C|\ge1$. Among the cobipartite supergraphs of $H$ on the same
vertex set, choose one with as many edges as possible subject to
\[
  \omega(H)\le k
  \qquad\text{and}\qquad
  C\text{ remaining attachable}.
\]
Such a supergraph exists, since the original graph is admissible.
It is still not $L$-colorable, since a coloring of it would color the
original graph. Keep the notation $H$ for this supergraph, and let
$F=\overline H$ be its bipartite nonedge graph. Choose a maximum
matching $M$ of $F$ minimizing
\[
  t=|M\cap E(F[A])|.
\]

For each $c\in C$, the set $N_F(c)$ is contained in $A$ and is
independent in $F$. Attachability of $C$ is equivalent to the
following \emph{cover condition}:
\begin{equation}\label{eq:cover}
  N_F(c)\text{ meets every edge of }F[A]
  \qquad(c\in C).
\end{equation}
Indeed, its complement in $A$ is $N_H(c)\setminus C$, which must
be a clique of $H$.

Deleting edges of $F$ corresponds to adding edges of $H$ and
preserves bipartiteness. If $M$ is retained, the matching number
cannot decrease, and it cannot increase under edge deletion. Thus
\eqref{eq:konig} shows that the clique number is unchanged. Our
choice of $H$ consequently forbids any deletion of an edge of $F$
that preserves both $M$ and the cover condition.

\begin{claim}\label{cl:FA}
The graph $F[A]$ consists of precisely the $t$ edges of $M$ with
both ends in $A$.
\end{claim}
\begin{proof}
Delete from $F[A]$ every edge not in $M$. The sets $N_F(c)$ are
unchanged, so the cover condition survives, and $M$ survives as
well. Edge maximality of $H$ implies that no edge was deleted.
\end{proof}

Write these $t$ edges as $x_jy_j$, $1\le j\le t$.

\begin{claim}\label{cl:CA}
For every $c\in C$ and $j$, exactly one of $cx_j,cy_j$ belongs to
$E(F)$. Apart from these covering edges, the only $C$--$A$ edges
of $F$ are the edges of $M$. Every $C$--$A$ edge of $M$ ends at a
vertex isolated in $F[A]$.
\end{claim}
\begin{proof}
The cover condition gives at least one of $cx_j,cy_j$. Both cannot
be present, since together with $x_jy_j$ they would form a triangle
in the bipartite graph $F$. Keep all these covering edges and all
edges of $M$, and delete any other $C$--$A$ edge. This preserves
the cover condition and $M$, so maximality permits no such edge.
Finally, an endpoint in $A$ of a $C$--$A$ matching edge cannot also
be an endpoint of one of the matching edges $x_jy_j$. The last
assertion follows from Claim~\ref{cl:FA}.
\end{proof}

\begin{claim}\label{cl:t}
We have $t\ge1$, and every vertex $c\in C$ is matched by $M$ to a
distinct vertex $z_c\in A$ isolated in $F[A]$.
\end{claim}
\begin{proof}
If $t=0$, Claim~\ref{cl:FA} says that $A$ is a clique of $H$.
Lemma~\ref{lem:cobi-asymmetric}, applied with first part $A$ and
second part $C$, gives an $L$-coloring, a contradiction. Hence
$t\ge1$.

Suppose that some $c\in C$ is unmatched. By
Claim~\ref{cl:CA}, it has an edge to one end of $x_jy_j$ for any
$j$. Replace $x_jy_j$ in $M$ by this edge. The new matching has
the same size as $M$, but has one fewer edge with both ends in $A$,
contradicting the choice of $M$. Thus every $c$ is matched. Its
partner lies in $A$ because $C$ is independent in $F$, and is
isolated in $F[A]$ by Claim~\ref{cl:CA}.
\end{proof}

\begin{claim}\label{cl:pattern}
The endpoints of the edges $x_jy_j$ can be labeled, and $C$ can be
partitioned into $C^+\dcup C^-$, so that in $F$ every vertex of
$C^+$ is adjacent to all $y_j$ and every vertex of $C^-$ is
adjacent to all $x_j$.
\end{claim}
\begin{proof}
Fix a bipartition $(R_0,R_1)$ of $F$, and label the endpoints so
that $x_j\in R_0$ and $y_j\in R_1$. Put
$C^+=C\cap R_0$ and $C^-=C\cap R_1$. For each $c$ and $j$,
Claim~\ref{cl:CA} supplies exactly one edge, and bipartiteness
forces its other endpoint to lie in the opposite part. This is
the claimed pattern.
\end{proof}

Label $C=\{c_1,\ldots,c_m\}$ and write $z_i=z_{c_i}$. Every
remaining vertex of $A$, outside the sets $\{z_i\}$ and
$\{x_j,y_j\}$, is isolated in $F$ by
Claims~\ref{cl:FA} and~\ref{cl:CA}; let these vertices be
$U=\{v_1,\ldots,v_u\}$. The preceding claims show that the
nonedges of $H$ are exactly those in \eqref{eq:J-nonedges}.
Thus $H=J(m,t,u;C^+,C^-)$ and
\[
  q=\omega(H)=m+t+u\le k.
\]
The lists on $C$ have at least $m$ colors, those on
$Z\cup X\cup Y$ have at least $k\ge q$, and those on $v_\ell$
have at least $k\ge q\ge m+\ell$. Hence
Lemma~\ref{lem:canonical-coefficient} gives an $L$-coloring of
$H$, the required contradiction.
\end{proof}

\begin{corollary}\label{cor:cobipartite-Mk}
Every cobipartite graph $H$ satisfies \Mk\ for every integer
$k\ge\omega(H)$.
\end{corollary}

\begin{proof}
In an \Mk-instance $(H,Q,r,L)$, the bound $|Q|+r\le k$ gives
$|L(q)|\ge k-r\ge |Q|$ on $Q$. All other lists have size at
least $k$. Apply Theorem~\ref{thm:M} with $C=Q$. This also covers
$Q=\varnothing$ and $r=0$.
\end{proof}

\section{Line graphs of bipartite multigraphs}\label{sec:linegraphs}

We prove the mixed-list property \Mk\ for line graphs of bipartite
multigraphs.  Galvin's kernel argument first gives a stronger list bound on
a full star.  A classification of attachable cliques then reduces the
general case to this full-star statement.

Throughout this section, $B$ is a bipartite multigraph with bipartition
$(P,P')$, and $G=\mathcal{L}(B)$.  Write $\st_B(w)$ for the set of edges incident with
$w$, counting parallel edges separately.  Every nonempty clique of $\mathcal L(B)$
is contained in a star: two nonparallel incident edges force a common
endpoint for all edges of the clique, since $B$ has no triangle; a family
of parallel edges already has a common endpoint.  Consequently
$\omega(\mathcal L(B))=\Delta(B)$, with both parameters taken to be zero when
$E(B)=\varnothing$.  Write $E(w,u)$ for the parallel edges joining $w$ and
$u$.  We call $u$ \emph{pendant relative to $w$} when $w$ is its only
distinct neighbor; this condition does not exclude parallel edges.

\subsection{Kernels and Galvin's orientation}

An \emph{orientation} of a graph $G$ replaces each edge by one or both of
its arcs.  A \emph{kernel} of a digraph is an independent set $K$ such that
every vertex outside $K$ has an arc into $K$; a digraph is
\emph{kernel-perfect} if every induced subdigraph has a kernel.  Outdegrees
$d^+(v)$ count distinct out-neighbors. We begin with the kernel lemma
of Bondy, Boppana and Siegel (see Alon and
Tarsi~\cite[Remark~2.4]{AlonTarsi1992}, or Bondy and
Murty~\cite[Theorem~14.19]{bondy-murty}).

\begin{lemma}[Kernel lemma]\label{lem:kernel}
Let $D$ be a kernel-perfect orientation of $G$.  If
$|L(v)|\ge d^+_D(v)+1$ for every $v$, then $G$ is $L$-colorable.
\end{lemma}

{The following lemma is essentially a retelling of \cite[Theorem 4.1]{Galvin}.}

\begin{lemma}[Galvin's orientation]\label{lem:galvin}
Let $c:E(B)\to\{1,\dots,k\}$ be a proper edge coloring.  Orient $\mathcal L(B)$ by
directing, at each vertex of $P$, from smaller to larger colors, and at each
vertex of $P'$, from larger to smaller colors (parallel edges get both arcs).
The result is a kernel-perfect orientation of $\mathcal L(B)$ with
$d^+(e)\le k-1$ for every edge $e$.  If $w\in P$ and $\st(w)$ receives exactly
the colors $1,\dots,q$, then $d^+(e)\le q-1$ for $e\in\st(w)$.
\end{lemma}

K\H{o}nig's edge-coloring theorem provides a proper $k$-edge-coloring whenever
$k\ge\Delta(B)$, and a permutation of the colors makes any prescribed star of
size $q$ receive $\{1,\dots,q\}$.

\begin{lemma}[Full stars]\label{lem:star}
Let $k\ge\Delta(B)$, and let either $Q=\varnothing$ or
$Q=\st_B(w)$ for some $w\in V(B)$.  Let $r\ge0$ satisfy $|Q|+r\le k$.
If $|L(e)|\ge k-r$ on $Q$ and $|L(e)|\ge k$ elsewhere, then $\mathcal{L}(B)$ is
$L$-colorable.
\end{lemma}

\begin{proof}
If $Q=\varnothing$, use any proper $k$-edge-coloring of $B$ and apply
Lemmas~\ref{lem:galvin} and~\ref{lem:kernel}; this is Galvin's ordinary
list-edge-coloring theorem.  Otherwise exchange the bipartition classes
if necessary so that $w\in P$, and permute the colors of a proper
$k$-edge-coloring so that $\st_B(w)$ receives $1,\dots,|Q|$.
Lemma~\ref{lem:galvin} gives outdegree at most $|Q|-1\le k-r-1$ on $Q$
and at most $k-1$ elsewhere. An application of the kernel lemma finishes the proof.
\end{proof}

\subsection{Attachable cliques of a bipartite line graph}

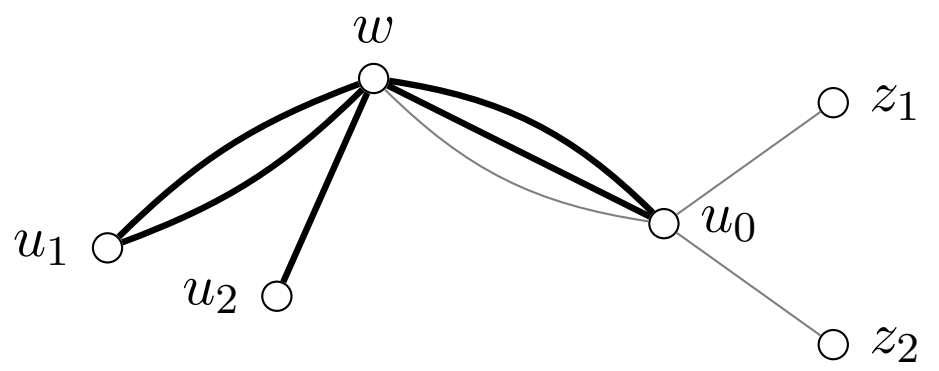
\begin{figure}[!htbp]
\centering
\begin{tikzpicture}[scale=0.8]
\node[vtx,label=above:$w$] (w) at (0,0) {};
\node[vtx,label=left:$u_1$] (u1) at (-2.2,-1.4) {};
\node[vtx,label=left:$u_2$] (u2) at (-0.8,-1.8) {};
\node[vtx,label=right:$u_0$] (u0) at (2.4,-1.2) {};
\node[vtx,label=right:$z_1$] (z1) at (3.8,-0.2) {};
\node[vtx,label=right:$z_2$] (z2) at (3.8,-2.2) {};
\draw[very thick] (w) to[bend left=12] (u1);
\draw[very thick] (w) to[bend right=12] (u1);
\draw[very thick] (w)--(u2);
\draw[very thick] (w) to[bend left=18] (u0);
\draw[very thick] (w) -- (u0);
\draw[gray] (w) to[bend right=18] (u0);
\draw[gray] (u0)--(z1); \draw[gray] (u0)--(z2);
\node[lab,align=center,text width=10cm] at (0.9,-3.25)
{Thick edges form $Q$; $u_1,u_2$ are pendant relative to $w$.\\
The bundle $E(w,u_0)$ has three edges, two in $Q_0$.\\
$N_0=\st(u_0)\setminus E(w,u_0)=\{u_0z_1,u_0z_2\}$.};
\end{tikzpicture}
\caption{An attachable clique in $\mathcal{L}(B)$, drawn in the root $B$.  Any edge at
$u_0$ other than $wu_0$ is nonadjacent (in $\mathcal{L}(B)$) to the edges $wu_1,wu_2$,
so \eqref{eq:alternative} forces everything outside $E(w,u_0)$ into $Q$.}
\label{fig:star}
\end{figure}

\begin{lemma}[Classification]\label{lem:star-classification}
Let $Q\ne\varnothing$ be an attachable clique of $\mathcal{L}(B)$.  Then $Q\subseteq\st(w)$
for some vertex $w$, and for every edge $wu\in Q$,
\begin{equation}\label{eq:alternative}
\st(w)\setminus\bigl(Q\cup E(w,u)\bigr)=\varnothing
\quad\text{or}\quad u\text{ is pendant relative to }w .
\end{equation}
Consequently at least one of the following holds:
\begin{enumerate}[label=\textup{(\roman*)}]
\item $Q=\st(w)$;
\item every far end $u$ of an edge $wu\in Q$ is pendant relative to $w$;
\item $w$ has exactly one nonpendant neighbor $u_0$, and
$Q=(\st(w)\setminus E(w,u_0))\cup Q_0$ with
$\varnothing\ne Q_0\subsetneq E(w,u_0)$ (Figure~\ref{fig:star}).
\end{enumerate}
\end{lemma}

\begin{proof}
Pairwise incident edges of a triangle-free multigraph share a vertex; this
gives $w$.  Let $wu\in Q$, let $wt\in\st(w)\setminus Q$ with $t\ne u$, and
let $ut'$ be an edge with $t'\ne w$.  Bipartiteness puts $t$ and $t'$ in
different parts, so $wt$ and $ut'$ are disjoint edges; they are neighbors of
$wu$ in $\mathcal{L}(B)$, nonadjacent, and outside $Q$ (as $ut'\notin\st(w)$).  This
contradicts attachability and proves \eqref{eq:alternative}.

Assume $Q\ne\st(w)$ and some far end $u_0$ of an edge of $Q$ is nonpendant.
Then \eqref{eq:alternative} for $wu_0$ gives $\st(w)\setminus E(w,u_0)\subseteq Q$,
so $Q_0=Q\cap E(w,u_0)$ is nonempty, and it is proper because $Q\ne\st(w)$.
Any other neighbor $u$ of $w$ has an edge in $Q$ (all of $\st(w)\setminus
E(w,u_0)$ lies in $Q$), and $E(w,u_0)\setminus Q_0\ne\varnothing$ lies
outside $Q\cup E(w,u)$, so \eqref{eq:alternative} makes $u$ pendant.  If no
such $u_0$ exists, (ii) holds.
\end{proof}

\subsection{The mixed-list lemma for bipartite line graphs}

\begin{lemma}[Mixed lists on bipartite line graphs]\label{lem:LG}
Let $G=\mathcal{L}(B)$ with $B$ bipartite, and $k\ge\omega(G)$.  Then $G$ satisfies
\Mk.
\end{lemma}

\begin{proof}
Let $(G,Q,r,L)$ be an \Mk-instance.  Note $\Delta(B)=\omega(G)\le k$.  If
$Q=\varnothing$ or $r=0$, all lists have at least $k$ colors and
Lemma~\ref{lem:star} with its explicit empty-set alternative applies.  Otherwise use
Lemma~\ref{lem:star-classification}.

\emph{Case (i).}  This is Lemma~\ref{lem:star}.

\emph{Case (ii).}  Color the clique $Q$ injectively from its lists, which have
at least $k-r\ge|Q|$ colors.  Since the far ends of $Q$ are pendant, the
vertices of $G-Q$ adjacent to $Q$ are exactly $Q'=\st(w)\setminus Q$.  Delete
the colors used on $Q$ from the lists on $Q'$; these lists keep at least
$k-|Q|$ colors.  In $B-Q$ (delete the edges of $Q$) the set $Q'$ is the full
star at $w$, and $|Q'|+|Q|=\deg_B(w)\le k$, so Lemma~\ref{lem:star} colors
$\mathcal{L}(B-Q)=G-Q$ compatibly.

\emph{Case (iii).}  Put $E_0=E(w,u_0)$ and $N_0=\st(u_0)\setminus E_0\ne\varnothing$.
First color the whole clique $\st(w)=Q\cup(E_0\setminus Q_0)$ from $L$.  Hall's condition \eqref{eq:hall} holds: a nonempty set
$T\subseteq Q$ has
$|L(T)|\ge k-r\ge|Q|\ge|T|$, and a set $T$ meeting $E_0\setminus Q_0$ has
$|L(T)|\ge k\ge\deg_B(w)\ge|T|$.  Let $f$ be this coloring.  The remaining
graph is $G-\st(w)=\mathcal{L}(B-w)$.  All neighbors of $w$ other than $u_0$ are
pendant, so the only vertices of $\mathcal{L}(B-w)$ adjacent to $\st(w)$ are those of
$N_0$, which are adjacent exactly to $E_0$.  Delete $f(E_0)$ from the lists on
$N_0$: they keep at least $k-|E_0|$ colors, $N_0$ is the full star at $u_0$ in
$B-w$, and $|N_0|+|E_0|=\deg_B(u_0)\le k$.  Lemma~\ref{lem:star} colors
$\mathcal{L}(B-w)$, and together with $f$ this is an $L$-coloring of $G$.
\end{proof}

\section{Elementary graphs}\label{sec:elementary}
We now prove \Pk\ for elementary graphs.  
In this section, we strengthen the chromatic-choosability theorem of Narayanan, Vasudevan and Somasundaram~\cite{NNS} to a 
precoloring extension form. We use their Deficiency Lemma directly. We then prove what we call the homogeneous-pair lemma.

\subsection{Representation and auxiliary results}

An edge is \emph{flat} if it lies in no triangle. \emph{Augmenting} a flat edge $xy$ means replacing $x$ and $y$ by disjoint nonempty cliques $X$ and $Y$, joined by at least one $X$--$Y$ edge, with $X$ complete to $N(x)\setminus\{y\}$ and $Y$ complete to $N(y)\setminus\{x\}$. The pair $(X,Y)$ is the \emph{augment}, and $X,Y$ are its \emph{shores}. Several flat edges forming a matching may be augmented simultaneously~\cite{MaffrayReed}. Adjacencies between different augments are inherited from the original vertices.

\begin{theorem}[Maffray--Reed~\cite{MaffrayReed}]\label{thm:MR}
A graph is elementary if and only if it is obtained from the line graph
of a bipartite multigraph by augmenting pairwise vertex-disjoint flat
edges. Elementary graphs are claw-free and perfect, and the class is
closed under induced subgraphs.
\end{theorem}

We call this a \emph{representation}, and its bipartite multigraph the
\emph{root}. By~\cite[Lemma~3.3]{NNS}, we may normalize the root so that
each augmented pair consists of edges $pw,wq$ with $p\ne q$ and
$\deg(w)=2$: an isolated pair of parallel edges is replaced by a
two-edge path.






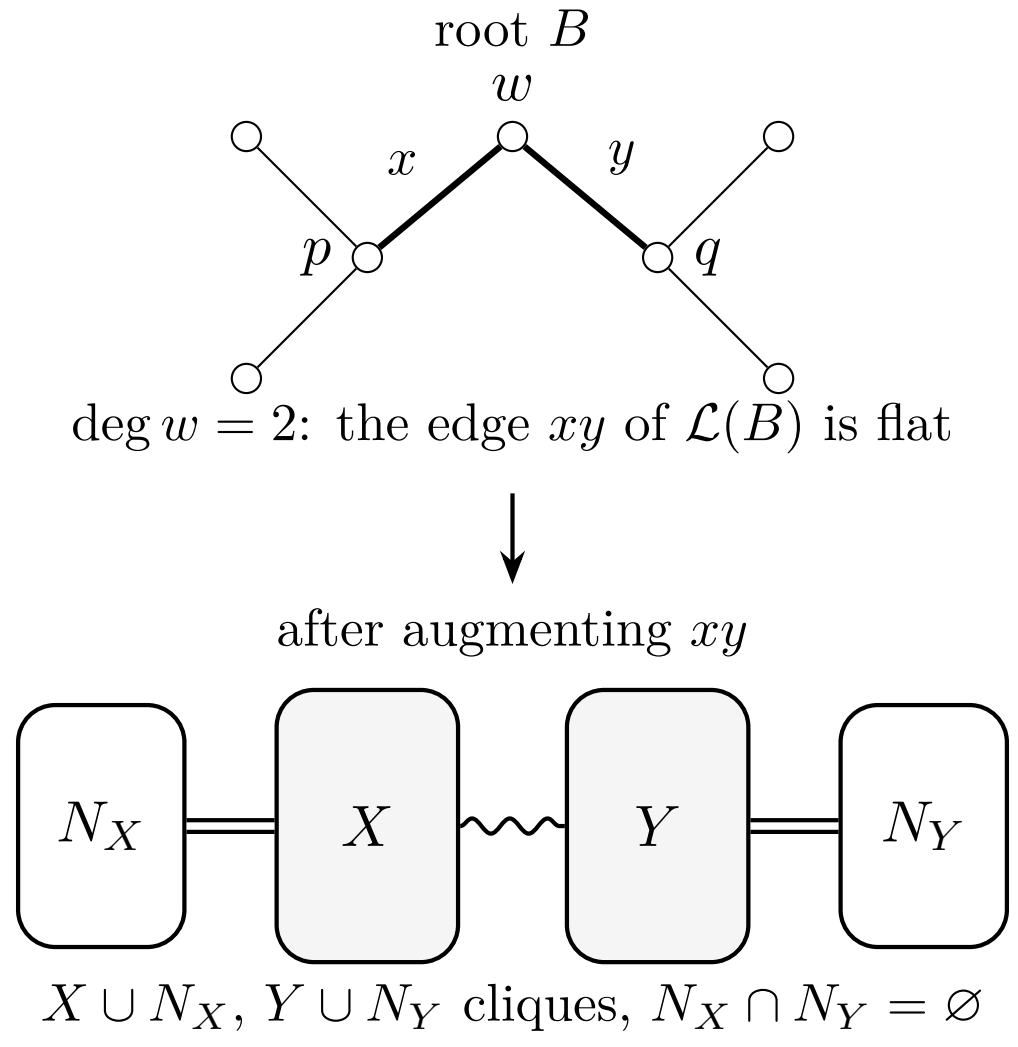
\begin{figure}[!htbp]
\centering
\begin{tikzpicture}[scale=.8]
\begin{scope}[xshift=0cm]
\node[lab] at (1.2,1.9) {root $B$};
\node[vtx,label=left:$p$] (p) at (0,0) {};
\node[vtx,label=above:$w$] (w) at (1.2,1.0) {};
\node[vtx,label=right:$q$] (q) at (2.4,0) {};
\node[vtx] (p1) at (-1.0,1.0) {}; \node[vtx] (p2) at (-1.0,-1.0) {};
\node[vtx] (q1) at (3.4,1.0) {}; \node[vtx] (q2) at (3.4,-1.0) {};
\draw[very thick] (p)--node[lab,above left]{$x$}(w)--node[lab,above right]{$y$}(q);
\draw (p)--(p1); \draw (p)--(p2); \draw (q)--(q1); \draw (q)--(q2);
\node[lab] at (1.2,-1.4) {$\deg w=2$: the edge $xy$ of $\mathcal L(B)$ is flat};
\end{scope}
\draw[arr] (1.2,-1.95)--(1.2,-2.7);
\begin{scope}[xshift=-1cm,yshift=-4.7cm]
\node[lab] at (2.2,1.6) {after augmenting $xy$};
\node[clq,fill=gray!8,minimum width=12mm,minimum height=18mm] (X) at (1.0,0) {$X$};
\node[clq,fill=gray!8,minimum width=12mm,minimum height=18mm] (Y) at (3.4,0) {$Y$};
\node[clq,minimum width=11mm,minimum height=16mm] (NX) at (-1.2,0) {$N_X$};
\node[clq,minimum width=11mm,minimum height=16mm] (NY) at (5.6,0) {$N_Y$};
\draw[full] (NX)--(X); \draw[full] (Y)--(NY);
\draw[incomp] (X)--(Y);
\node[lab] at (2.2,-1.5) {$X\cup N_X$, $Y\cup N_Y$ cliques, $N_X\cap N_Y=\varnothing$};
\end{scope}
\end{tikzpicture}
\caption{Augmenting a flat edge.  Double lines denote complete pairs, the wavy
line a specified $X$--$Y$ adjacency with at least one cross edge.
The graph on $X\cup Y$ is connected; the bipartite graph of cross edges
need not be connected.}
\label{fig:augment}
\end{figure}

\Needspace{17\baselineskip}
\begin{lemma}[Augments]\label{lem:augment}
Fix a normalized representation of an elementary graph $G$.
For an augment $(X,Y)$, put $N_X=N(X)\setminus Y$ and
$N_Y=N(Y)\setminus X$.
\begin{enumerate}[label=\textup{(\alph*)}]
\item Every vertex outside $X\cup Y$ is complete or anticomplete to each of
$X$ and $Y$; $X\cup N_X$ and $Y\cup N_Y$ are cliques; and $N_X\cap N_Y=\varnothing$,
so no outside vertex sees both $X$ and $Y$.
\item Adding edges between $X$ and $Y$ preserves  elementarity. 
\item Any collection of complete augments can be absorbed into the root,
leaving the remaining augments unchanged. Moreover, if $G'$ is an induced
subgraph in which the surviving sides of every augment form a complete
pair, then $G'$ is a bipartite line graph; surviving sides may be empty.
\end{enumerate}
\end{lemma}

\begin{proof}
Part (a) is~\cite[Lemma~3.3]{NNS}, and (b) follows from the representation.
For (c), use the construction of~\cite[Lemma~3.8]{NNS}: replace
$pw,wq$ by parallel edges of the required sizes, allowing size zero for a deleted side. Delete root edges corresponding
to deleted unaugmented vertices. Note that, this process can be done partially as well. That is, replace only the complete augments to be absorbed.
 The center of another designated flat pair cannot be $p$ or $q$ as  its
 degree would be two, so its support would contain one of $pw,wq$,
contrary to disjointness of the supports.  Thus the other designated
 pairs remain flat and keep their original augments.  This proves the
 partial-absorption assertion as well.

\end{proof}

We need the deficiency lemma proved in~\cite{NNS}, which we restate here adjusting the notation.

\begin{lemma}[{\cite[Lemma~2.5]{NNS}}]\label{lem:deficiency}
Let $H=K_{a+b}-M$, where the vertex set is $X\dcup Y$ with $|X|=a$, $|Y|=b$,
and $M$ is a matching of $\mu$ pairs $xy$ with $x\in X$, $y\in Y$.  Let
$X_0,Y_0$ be the sets of matched vertices and $X_*=X\setminus X_0$,
$Y_*=Y\setminus Y_0$. Let $R$ be lists with
$|R(x)|\ge a$ on $X$ and $|R(y)|\ge b$ on $Y$.  If $H$ is not $R$-colorable,
there is a set $T\subseteq X_*\cup Y_*$ meeting both $X_*$ and $Y_*$ with
\begin{equation}\label{eq:deficiency}
|R(T)|\le\mu+|T|-1 .
\end{equation}
In particular, $H$ is $R$-colorable if either unmatched part is empty.
\end{lemma}

\subsection{A homogeneous-pair obstruction}

   Fix $k\ge0$ and a hereditary class $\mathcal C$ of claw-free perfect graphs. An \emph{extremal counterexample} $(G,S,\phi,L)$ is chosen as follows:

\begin{enumerate}[label=\textup{(\roman*)}]
\item minimize $|V(G)|$ over all uncolorable
\Pk-instances with $G\in\mathcal C$ and $\omega(G)\le k$;
\item keeping $V(G)$, $S$, $\phi$ and $L$ fixed, maximize $|E(G)|$ among all uncolorable instances on this vertex set with $G\in\mathcal C$, $\omega(G)\le k$, and $S$ a clique of simplicial vertices.
\end{enumerate}
By Lemma \ref{lem:minimal}? every proper induced subgraph of $G$ satisfies \Pk.


\begin{lemma}[Homogeneous-pair lemma]\label{lem:GP9}
Let $(G,S,\phi,L)$ be an extremal counterexample chosen as above. There are no disjoint nonempty cliques $A_1,A_2\subseteq V(G)\setminus S$
such that:
\begin{enumerate}[label=\textup{(\roman*)}]
\item every vertex outside $A_1\cup A_2$ is complete or anticomplete
to each $A_i$;
\item writing $Q_i=N(A_i)\setminus(A_1\cup A_2)$, each $A_i\cup Q_i$
is a clique;
\item $G[A_1\cup A_2]$ is not complete, and $G+uv\in\mathcal C$
for every missing pair $u\in A_1$, $v\in A_2$.
\end{enumerate}
\end{lemma}

\begin{proof}
We adapt the list-copying argument of~\cite[Proposition~3.7]{NNS}.
Put $a_i=|A_i|$ and $W=Q_1\cap Q_2$; then $|Q_i|\le k-a_i$.

\begin{figure}[!htbp]
\centering
\begin{tikzpicture}[scale=0.9]
\node[cbox,minimum width=16mm,minimum height=9mm] (SU) at (0,1.9) {$S_U$};
\node[cbox,minimum width=16mm,minimum height=9mm] (M1) at (0,0.6) {$M_1$};
\node[cbox,minimum width=16mm,minimum height=9mm] (A1s) at (0,-0.7) {$A_1^*$};
\node[clq,fit=(SU)(M1)(A1s),label=below:{$U=S_U\cup A_1$}] (U) {};
\node[cbox,minimum width=16mm,minimum height=9mm] (M2) at (5.0,0.6) {$M_2$};
\node[cbox,minimum width=16mm,minimum height=9mm] (A2s) at (5.0,-0.7) {$A_2^*$};
\node[clq,fit=(M2)(A2s),label=below:{$V=A_2$}] (V) {};
\draw[nonedge] (M1.east)--(M2.west);
\node[lab] at (2.55,1.42) {$M$ (missing pairs)};
\draw[full] (A1s.east)--(A2s.west);
\draw[full] (A1s.east)--(M2.west);
\draw[full] (M1.east)--(A2s.west);
\node[cbox,minimum width=20mm,minimum height=9mm] (NU) at (-0.2,3.8) {$N_U\setminus W$};
\node[cbox,minimum width=14mm,minimum height=9mm] (W) at (2.5,3.8) {$W$};
\node[cbox,minimum width=20mm,minimum height=9mm] (NV) at (5.2,3.8) {$N_V\setminus W$};
\draw[full] (U.north)--(NU.south);
\draw[full] (V.north)--(NV.south);
\draw[full] (W.south west)--(U.north east);
\draw[full] (W.south east)--(V.north west);
\node[lab,align=center,text width=11.5cm] at (2.5,-2.65)
{$Q_1=N_U\cup S_U$, \quad $Q_2=N_V$; \quad $A_i\cup Q_i$ are cliques.\\
$M_i$ are the matched vertices and $A_i^*$ the unmatched vertices.\\
$K^*=M_1\cup A_1^*\cup A_2^*\cup W$ is a clique of $G-M_2$.};
\end{tikzpicture}
\caption{The homogeneous pair of Lemma~\ref{lem:GP9} after saturation.  Double
lines are complete pairs; $S_U$ is anticomplete to $V$; the only missing pairs
between $A_1$ and $A_2$ form the matching $M$.}
\label{fig:GP9}
\end{figure}

\emph{Saturation.}
Let $F$ be the bipartite graph of missing $A_1$--$A_2$ pairs, and let
$M$ be a maximum matching of $F$, of size $\mu\ge1$.
If $uv\in E(F)\setminus M$, adding $uv$ leaves the clique number of
$G[A_1\cup A_2]$ unchanged by Lemma~\ref{lem:saturation}.
Every new clique lies in $A_1\cup A_2\cup W$. Since $W$ is a clique complete to $A_1\cup A_2$, a maximum clique of $G[A_1\cup A_2]$ together with $W$ is a clique of $G$, so $\omega(G[A_1\cup A_2])+|W|\le\omega(G)\le k$; hence $\omega(G+uv)\le k$.

The addition preserves simpliciality and the prescribed coloring of
$S$, since its ends lie outside $S$. Together with (iii), this contradicts
edge maximality. Hence $F=M$.

Let $M_i$ be the matched vertices of $A_i$, and put
$A_i^*=A_i\setminus M_i$, $\alpha_i=|A_i^*|$, so $a_i=\mu+\alpha_i$.
Both $A_i^*$ are nonempty: otherwise color $G-(A_1\cup A_2)$ by
minimality, extending $\phi$, and remove the colors on $Q_i$ from lists
on $A_i$. The remaining lists have size at least $a_i$, and
Lemma~\ref{lem:deficiency} completes the coloring.

\emph{Copying a list.}
Set
\[
m_i=\min_{\varnothing\ne T\subseteq A_i^*}(|L(T)|-|T|),
\]
and interchange the indexed data so that
$m_1+\alpha_1\le m_2+\alpha_2$.
Choose $\varnothing\ne Z\subseteq A_1^*$ attaining $m_1$, and $z\in Z$.
On $G'=G-M_2$, replace every list on $M_1$ by $L(z)$.
Minimum order gives a coloring $f$ extending $\phi$: $S$ is retained
and all unprescribed lists still have size $k$.
Put $F_i=f(Q_i)$ and $R(v)=L(v)\setminus F_i$ for $v\in A_i$.
Then $|R(v)|\ge a_i$. An $R$-coloring of $G[A_1\cup A_2]$ would
combine with $f$ outside the pair to color $G$. Thus
Lemma~\ref{lem:deficiency} gives $T\subseteq A_1^*\cup A_2^*$ meeting
both parts with $|R(T)|\le\mu+|T|-1$.

Write $t_i=|T\cap A_i^*|$. Since
$L(T\cap A_2^*)\subseteq R(T)\cup F_2$,
\[
m_2+t_2\le\mu+t_1+t_2-1+k-a_2,
\qquad t_1\ge m_2+\alpha_2-k+1.
\]
But $|L(Z)|=m_1+|Z|\ge k$, so $T\cap Z=\varnothing$ would give
\[
t_1\le\alpha_1-|Z|\le m_1+\alpha_1-k
\le m_2+\alpha_2-k,
\]
a contradiction. Put $I=T\cap Z\ne\varnothing$ and $T'=T\cup Z$.
The choice of $Z$ gives
$|L(Z)\setminus L(I)|\le|Z\setminus I|$.
All lists on $Z$ lose the same set $F_1$, whence
\[
|R(T')|\le |R(T)|+|Z\setminus T|\le\mu+|T'|-1.
\]
On the other hand, $M_1\cup T'$ is a clique of $G'$.
Its vertices in $M_1$ receive distinct colors from
$L(z)\setminus F_1=R(z)\subseteq R(T')$, and those in $T'$ receive
colors from their unchanged residual lists. Thus $f$ assigns
$\mu+|T'|$ distinct colors from $R(T')$, a contradiction.
\end{proof}

\subsection{Precoloring extension}

\begin{theorem}[Elementary precoloring extension]\label{thm:elementary-P}
Every elementary graph $G$ satisfies \Pk\ for every $k\ge\omega(G)$.
\end{theorem}

\begin{proof}
Fix $k$ and suppose the assertion fails. Choose an extremal
counterexample in the elementary class.
By Lemma~\ref{lem:minimal}, $G$ is connected. Set $Q=N(S)$, with
$Q=\varnothing$ when $S=\varnothing$.
Lemma~\ref{lem:simplicial-clique} reduces the extension problem to an
\Mk-instance $(G-S,Q,|S|,L')$. In particular, $G$ is not cobipartite,
by Corollary~\ref{cor:cobipartite-Mk}.

Fix a normalized representation. The clique $S$ meets at most one side
of each augment $(X,Y)$: vertices of $S$ are true twins, whereas a
vertex of $N_X$ sees $X$ and misses $Y$, and conversely for $N_Y$.
If $S$ met both sides, therefore, $N_X=N_Y=\varnothing$; connectedness
would give $G=G[X\cup Y]$, a contradiction. Orient each augment so that
$S\cap Y=\varnothing$, and put $S_X=S\cap X$.

Suppose that $G[(X\setminus S_X)\cup Y]$ is not complete.
If $S_X=\varnothing$ or $N_X\ne\varnothing$, then $S_X$ is
anticomplete to $Y$: a neighbor in $Y$ and a vertex of $N_X$ would be
nonadjacent neighbors of a simplicial vertex in $S_X$. In particular, $X\setminus S_X \neq \varnothing$, because the augment has at least one $X--Y$ edge and none of these edges leaves $S_X$.
The pair $(A_1,A_2)=(X\setminus S_X,Y)$ now contradicts
Lemma~\ref{lem:GP9}, with $Q_1=N_X\cup S_X$ and $Q_2=N_Y$;
edge additions are allowed by Lemma~\ref{lem:augment}(b).

It remains that $S_X\ne\varnothing$ and $N_X=\varnothing$.
Let $Y_S=N(s)\cap Y$ for $s\in S_X$; this is independent of $s$
because the vertices of $S$ are twins. Simpliciality makes $Y_S$
complete to $X$. Hence every missing pair of $G[(X\setminus S_X)\cup Y]$ joins $A_1=X\setminus S_X$ to $A_2=Y\setminus Y_S$, and both sets are nonempty because at least one such pair exists.

These are homogeneous cliques with
\[
Q_1=S_X\cup Y_S,\qquad Q_2=Y_S\cup N_Y.
\]
Both $A_i\cup Q_i$ are cliques, and their missing edges lie within
the original augment. Lemma~\ref{lem:GP9} again gives a contradiction.

Consequently every augment becomes complete after deleting $S$.
By Lemma~\ref{lem:augment}(c), $G-S$ is a bipartite line graph.
Lemma~\ref{lem:LG} colors the reduced \Mk-instance, and restoring
$\phi$ gives the required extension.
\end{proof}

\begin{remark}
The attachment result in Section~\ref{sec:elementary-atom-attachment}
will establish when adjoining a prescribed simplicial clique preserves
elementarity; this closure property is not asserted here.
\end{remark}

\section{Graphs of interval type}\label{sec:interval}

The interval-type argument uses the ordered clique structure of the
graphs of interval type introduced by Chudnovsky--Seymour~\cite{CS}.
After saturating each incomplete pair, we extract a coefficient one
fiber at a time.  Processing the first end of an incomplete pair
transfers a Vandermonde factor to its other end; processing that end
replaces the factor by a nonzero factorial.  This local transfer is the
reason crossing incomplete pairs can be handled.

For variables or constants $\zeta_1,\dots,\zeta_h$, write
\[
\Van(\zeta_1,\dots,\zeta_h)=\prod_{i<j}(\zeta_j-\zeta_i)
=\det\bigl[\zeta_i^{\,j-1}\bigr]_{i,j=1}^h,\qquad \Van(\varnothing)=1 .
\]
We use vertex names also for their variables.

\Needspace{10\baselineskip}
\subsection{Fibers, windows and saturation}

\begin{definition}\label{def:interval}
A graph $G$ is \emph{of interval type} if its vertex set has an ordered
partition into nonempty cliques $C_1,\dots,C_m$, called \emph{fibers},
such that:
\begin{enumerate}[label=\textup{(I\arabic*)}]
\item\label{condition: interval} if $u\in C_i$, $v\in C_j$, $w\in C_\ell$ with $i<j<\ell$ and $uw\in E(G)$,
then $uv,vw\in E(G)$;
\item\label{condition: incomplete-partner} each fiber forms an incomplete pair with at most one other fiber,
its \emph{incomplete partner}.
\end{enumerate}
\end{definition}

Here an \emph{incomplete pair} is neither complete nor anticomplete.
We allow $m=0$ for the empty graph.  Intersecting the fibers with an
induced subgraph and discarding empty intersections preserves (I1) and
(I2); a former incomplete pair may become complete or anticomplete.
Reversing the order also preserves both conditions. Thus the class is
hereditary, and we may process a representation from either end.


\Needspace{14\baselineskip}
\begin{lemma}[Windows]\label{lem:windows}
Let $(U,V)=(C_i,C_j)$, $i<j$, be an incomplete pair and define
$W=C_{i+1}\cup\dots\cup C_{j-1}$ to be its \emph{window}.
\begin{enumerate}[label=\textup{(\alph*)}]
\item $W$ is a clique complete to $U\cup V$.
\item No fiber before $U$ has a neighbor in $V$ and no fiber after $V$ has
a neighbor in $U$.  Hence every clique meeting $U$ and $V$ lies in
$U\cup W\cup V$, and the neighbors of $U$ in later fibers lie in $W\cup V$.
\item If two fibers $A,B$ are such that $A,B$ is a complete pair, then every fiber in between together with $A,B$ forms a clique.
\end{enumerate}
\end{lemma}

\begin{proof}
(a) First, we show that every edge between $W$ and $U\cup V$ exists in $G$.
Since $U,V$ forms an incomplete pair, there exists $u\in U,v\in V$ such that 
$uv\in E(G)$. This forces $u,v$ to be adjacent to every vertex in $W$, by 
\ref{condition: interval}. Since $U,V$ are the unique incomplete partners of each other, every other fiber is a complete or anti-complete pair for $U,V$, by \ref{condition: incomplete-partner}. The existence of one edge rules out anti-completeness and so, every fiber in $W$ forms a complete pair with $U,V$. By \ref{condition: interval}, $W$ forms a clique and since every $C_r$ in $W$ is a complete pair of $U,V$, $W$ is a complete pair of $U\cup V$.

(b) If a fiber $D$ before $U$ had a neighbor in $V$, then $(D,V)$
would be complete, since $U$ is the unique incomplete partner of $V$.
An edge from $D$ to the endpoint in $V$ of a missing $U$--$V$ pair
would then force that missing edge by \ref{condition: interval}, a
contradiction. Reversing the order proves the assertion for fibers after
$V$. The remaining assertions follow immediately.

(c) WLOG, assume $A$ appears before $B=C_r$ such that $A,B$ is a complete pair. So $C_{r-1}$ forms a complete pair with $A$, by \ref{condition: interval}. Recursively, each fiber between $A$ and $C_r$ is a complete pair with both $A$ and $C_r$ and this proves the assertion.
\end{proof}


\begin{lemma}[Matching saturation]\label{lem:interval-saturation}
Let $G$ be of interval type.  Adding, for each incomplete pair $(U,V)$, all
missing $U$--$V$ pairs outside a maximum matching $M_{UV}$ of the missing pairs
yields a graph of interval type with the same ordered fibers and the
same clique number.  In the resulting graph, the missing pairs between
the ends of each incomplete pair form a nonempty matching.  If a vertex set $S$ meets no endpoint fiber of an incomplete pair,
no edge at $S$ is added.
\end{lemma}

\begin{proof}
Let $G'$ be the graph after adding the edges mentioned in the statement of the lemma.
It is clear that there are no new complete or anti-complete pairs. Suppose $i<j<k$ and let
$uw\in E(G')\setminus E(G)$ with $u\in C_i,w\in C_k$. This implies that $C_i,C_k$
is an incomplete pair and we may apply Lemma \ref{lem:windows}(a) to get that $uv,vw\in E(G)$ for any $v\in C_j$. This gives us Condition \ref{condition: interval} showing that $G'$ is interval type. 

Let $u\in C_i, v\in C_j$ with $uv\notin E(G)$. We claim that adding $uv$ does not produce a larger clique if the edge addition does not decrease the size of the maximum matching on the complement of $G[C_i\cup C_j]$. If such a larger clique $K$ is formed, it necessarily contains $u$ and $v$. From Lemma \ref{lem:windows}(b), it is clear that the common neighbors of $u,v$ are a subset of $C_i\cup W \cup C_j$ where $W$ is the window of $C_i,C_j$ and so, $K\subseteq C_i\cup W\cup C_j$. Since $W$ is complete to $C_i\cup C_j$ (by Lemma \ref{lem:windows}), the largest clique $K$ is entirely controlled by the clique number of $G[C_i\cup C_j]$. This contradicts Lemma \ref{lem:saturation}. Therefore, adding the edges one by one shows that the clique number does not change. The rest of the assertions are straightforward.
\end{proof}

A graph of interval type is \emph{saturated} if the missing edges of every
incomplete pair form a matching; we then say the pair $(U,V)$ has
\emph{matched} vertices $a_1,\dots,a_\mu\in U$ and $b_1,\dots,b_\mu\in V$,
with $a_ib_i\notin E(G)$, and \emph{unmatched} vertices otherwise.

\begin{figure}[ht]
\centering
\begin{tikzpicture}[scale=.90]
\node[vtx,label=left:$a_1$] (a1) at (0,1.2) {};
\node[vtx,label=left:$a_2$] (a2) at (0,0) {};
\node[vtx,label=left:$u_1$] (u1) at (0,-1.2) {};
\node[clq,fit=(a1)(a2)(u1),label=below:$A$] (A) {};
\node[vtx,label=right:$b_1$] (b1) at (6.7,1.2) {};
\node[vtx,label=right:$b_2$] (b2) at (6.7,0) {};
\node[vtx,label=right:$z_1$] (z1) at (6.7,-1.2) {};
\node[clq,fit=(b1)(b2)(z1),label=below:$B$] (B) {};
\node[cbox,minimum width=12mm,minimum height=9mm] (Wb) at (3.2,2.6) {$W$};
\draw[full] (A.north)--(Wb.south west); \draw[full] (B.north)--(Wb.south east);
\node[cbox,minimum height=20mm,minimum width=10mm] (T) at (9.4,0) {$T'$};
\node[lab,below=1mm of T] {$|T'|=R$};
\draw[full] (B)--(T);
\draw[nonedge] (a1)--(b1); \draw[nonedge] (a2)--(b2);
\node[lab] at (3.2,-0.6) {all other $A$--$B$ pairs are edges};
\draw[arr,gray!50!black] (0.8,-2.2) to[bend right=20] node[lab,below]{$\Van(b_1,b_2)$} (5.6,-2.2);
\node[lab,align=left] at (4.6,-3.9) {sending at $A$ ($n=|W|+|B|$, $\mu=2$, $\alpha=1$):\quad
  $d(a_1)=d(a_2)=n+\alpha-1$,\ \ $d(u_1)=n$\\
  receiving at $B$ ($\mu=2$, $\beta=1$):\quad
  $d(b_1)=d(b_2)=R+\mu+\beta-1=R+2$,\ \ $d(z_1)=R$};
\end{tikzpicture}
\caption{Eliminating the incomplete pair $(A,B)$ with window $W$.  Dashed:
the matched nonedges $a_1b_1$, $a_2b_2$.  The window is complete to $A\cup B$;
$T'$ is the set of neighbors of $B$ in later fibers.  Eliminating $A$ leaves
the factor $\Van(b_1,b_2)$, which is absorbed when $B$ is eliminated.}
\label{fig:transfer}
\end{figure}
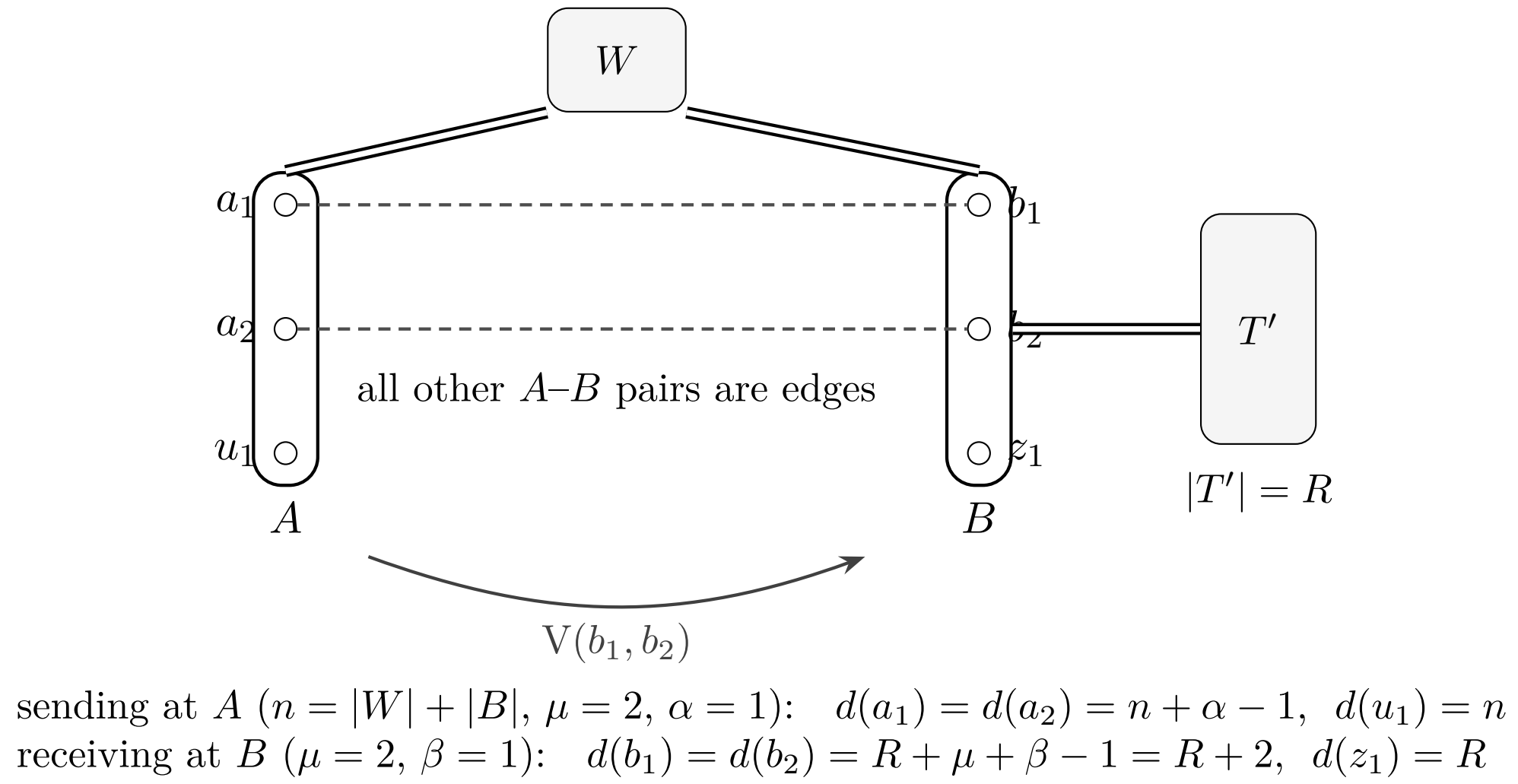

Consider a saturated pair $(A,B)$ with matched vertices $a_i\in A$, $b_i\in B$
($1\le i\le\mu$) and unmatched vertices $u_1,\dots,u_\alpha\in A$
(Figure~\ref{fig:transfer}).  Set $T=W\cup B$, where $W$ is the window of $(A,B)$, and put $n=|T|$.
Every neighbor of $A$ in a later fiber belongs to $T$.  The only
nonadjacent pairs between $A$ and $T$ are $a_i b_i$; in particular,
$T$ includes $b_i$ even if it has no neighbor in $A$.  Define
\[
P(x)=\prod_{\tau\in T}(x-\tau).
\]
Up to sign, the factors of the graph polynomial that involve $A$ and vertices
of $A\cup T$ are
\begin{equation}\label{eq:FA}
F_A=\Van(a_1,\dots,a_\mu,u_1,\dots,u_\alpha)\
\prod_{i=1}^{\mu}\frac{P(a_i)}{a_i-b_i}\
\prod_{r=1}^{\alpha}P(u_r),
\end{equation}
the quotients being polynomials since $b_i\in T$.


\Needspace{8\baselineskip}
\begin{lemma}[Sending identity]\label{lem:sending}
Let $d(a_i)=n+\alpha-1$ and $d(u_r)=n+r-1$.  Then
\[
\Bigl[\prod_i a_i^{d(a_i)}\prod_r u_r^{d(u_r)}\Bigr]F_A
=\pm\Van(b_1,\dots,b_\mu),
\]
the coefficient being taken in the variables of $A$.  All these exponents are
at most $\omega(G[A\cup T])-1$.  The identity holds verbatim if some elements
of $T\setminus\{b_1,\dots,b_\mu\}$ are constants.
\end{lemma}

\begin{proof}
Write $m=\mu+\alpha$ and $P(x)=x^n+c_1x^{n-1}+\cdots+c_n$ with $c_0=1.$
Each term of the Vandermonde determinant assigns a permutation of
$0,\dots,m-1$ to the variables. Since $Q_i(x)=P(x)/(x-b_i)$ has degree $n-1$, a term contributing to the
desired coefficient must assign Vandermonde exponent at least
$\alpha$ to each $a_i$. Thus the matched variables receive precisely
the exponents $\alpha,\dots,\alpha+\mu-1$, and the unmatched variables
receive $0,\dots,\alpha-1$.

Similarly, since $P$ has degree $n$, the exponent assigned to $u_r$
must be at least $r-1$. Working backwards from $r=\alpha$ forces this
exponent to be exactly $r-1$. Therefore each factor $P(u_r)$
contributes only its leading term $u_r^n$.

Let $Q$ be the $\mu\times\mu$ matrix with rows indexed by
$i=1,\dots,\mu$ and columns indexed by $j=0,\dots,\mu-1$, where
$Q_{i,j}=[x^{n-1-j}]Q_i(x)$. For each $\sigma\in S_\mu$, the corresponding Vandermonde term assigns
exponent $\alpha+\sigma(i)-1$ to $a_i$ and $Q_i(a_i)$ must contribute exponent $n-\sigma(i)$, whose
coefficient is $Q_{i,\sigma(i)-1}$.
Up to a common sign independent of $\sigma$, the contribution is
therefore
$\operatorname{sgn}(\sigma)\prod_{i=1}^{\mu}Q_{i,\sigma(i)-1}$.
Summing over all permutations, the desired coefficient is, up to sign,
\[
\sum_{\sigma\in S_\mu}
\operatorname{sgn}(\sigma)
\prod_{i=1}^{\mu}Q_{i,\sigma(i)-1}
=\det Q.
\]

Polynomial division gives $Q_{i,j}=\sum_{\ell=0}^{j}c_\ell b_i^{j-\ell}=b_i^j+c_1b_i^{j-1}+\cdots+c_j.$
Hence $Q=V_BC$, where $V_B$ is the Vandermonde matrix with entries
\(
(V_B)_{i,k}=b_i^k
\qquad
(1\le i\le\mu,\ 0\le k\le\mu-1),
\)
and
\[
C_{k,j}=
\begin{cases}
c_{j-k}, & k\le j,\\
0, & k>j,
\end{cases}
\qquad 0\le k,j\le\mu-1.
\]
Since $C$ is upper triangular with diagonal entries $c_0=1$, we obtain $\det Q=\det V_B.$

Finally, the saturated graph on $A\cup T$ is complete except for the
$\mu$ disjoint missing pairs $a_ib_i$. Its clique number is therefore
\[
|A|+n-\mu=n+\alpha.
\]
Every prescribed exponent is at most $n+\alpha-1$, with equality
attained. Thus the largest prescribed exponent is one less than the
clique number of this block.
\end{proof}

For the first nonconstant transfer, take $\mu=2$, $\alpha=0$ and
$T=\{b_1,b_2\}$.  The identity
reads $[a_1a_2]\,(a_2-a_1)(a_1-b_2)(a_2-b_1)=b_2-b_1$.  The output involves
only the matched vertices of $B$ and nothing in the window; this is what
makes crossing incomplete pairs independent.

Now let $B$ have matched vertices $b_1,\dots,b_\mu$ carrying the factor
$\Van(b_1,\dots,b_\mu)$, unmatched vertices $z_1,\dots,z_\beta$, and let $T'$
be the set of neighbors of $B$ in later fibers, with $R=|T'|$; by
Lemma~\ref{lem:windows}(c), $B\cup T'$ is a clique.  The factors involving
$B$, including the received one, are up to sign
\begin{equation}\label{eq:FB}
F_B=\Van(b_1,\dots,b_\mu,z_1,\dots,z_\beta)\,\Van(b_1,\dots,b_\mu)
\prod_{v\in B}\prod_{\tau\in T'}(v-\tau).
\end{equation}

\begin{lemma}[Receiving identity]\label{lem:receiving}
Let $d(b_i)=R+\mu+\beta-1$ and $d(z_r)=R+r-1$.  Then
\[
\Bigl[\prod_i b_i^{d(b_i)}\prod_r z_r^{d(z_r)}\Bigr]F_B=\pm\,\mu!\,,
\]
and all these exponents are at most $\omega(G[B\cup T'])-1$.  The case
$\mu=0$ (no received factor) is included, and elements of $T'$ may be
constants.
\end{lemma}

\begin{proof}
Set $m=\mu+\beta$. The exponent sum $\mu(R+\mu+\beta-1)+\beta R+\binom\beta2$ equals the degree
of $F_B$ in the variables of $B$, so only the top homogeneous part contributes.
Therefore each factor $v-\tau$ contributes $v$, whether $\tau$ is a variable or a
constant.  Removing $\prod_{v\in B}v^R$ leaves
\[
\Van(b)^2\,\Van(z)\prod_{i,r}(z_r-b_i),
\]
and we need the coefficient of $\prod_ib_i^{\mu+\beta-1}\prod_rz_r^{\,r-1}$. $\Van(z)$ contributes
monomials of total degree $\binom{\beta}{2}$, which coincides with the total degree of $\prod_r z_r^{r-1}$. Therefore, the cross difference factors should only
contribute $-b_i$. In the remaining $[\prod_ib_i^{\mu-1}]\Van(b)^2$, expanding both
determinants, once a permutation $\sigma$ of exponents is chosen in the first,
the second must use $\mu-1-\sigma(i)$ at $b_i$, i.e.\ the reversal of
$\sigma$; the $\mu!$ terms all have the same sign $(-1)^{\binom\mu2}$.  Hence
the coefficient is $\pm\mu!$.
Every exponent is nonnegative and at most $R+|B|-1$.
Since $B\cup T'$ is a clique, this upper bound is exactly
$\omega(G[B\cup T'])-1$.  For $\mu=0$ the matched block and its
Vandermonde factor are empty, and the coefficient is $0!=1$.
\end{proof}

\subsection{The elimination lemma}

Let $G$ be saturated of interval type with fibers $C_1,\dots,C_m$, and let
$p\in\{1,\dots,m\}$.  The \emph{elimination order toward $C_p$} is defined to be the
sequence
\[
C_1,\ C_2,\ \dots,\ C_{p-1},\ C_m,\ C_{m-1},\ \dots,\ C_{p+1};
\]
the \emph{full elimination order} is defined to be $C_1,\dots,C_m$.
 \begin{lemma}[Elimination]\label{lem:elimination}
Let $G$ be saturated of interval type with fibers $C_1,\dots,C_m$ and let
$\omega=\omega(G)$.  Let either $Z=\varnothing$, or $Z=C_p$ for a fiber
$C_p$ that has no incomplete partner and is straddled by no incomplete pair
(no pair $(C_i,C_j)$ with $i<p<j$).  Let $\phi:Z\to\RR$ be any map, and form $P_{G,\phi}$ by the product
formula in Corollary~\textup{\ref{cor:nullstellensatz}}, omitting edges
internal to $Z$.  Properness of $\phi$ is not needed for the following
algebraic identity; it will be required when the identity is used to
extend a coloring.
Then there are integers $0\le d(v)\le\omega-1$, $v\notin Z$, with
$\sum_{v\notin Z}d(v)=\deg P_{G,\phi}$ and
\[
\Bigl[\prod_{v\notin Z}x_v^{d(v)}\Bigr]P_{G,\phi}
=\pm\prod_{(U,V)}\mu(U,V)!\ \neq0,
\]
the product over all incomplete pairs, $\mu(U,V)$ being the size of the
matching of missing edges.
\end{lemma}

\begin{proof}
Our strategy will be as follows. We will compute
$\bigl[\prod_{v\notin Z}x_v^{d(v)}\bigr]P_{G,\phi}$ fiber by fiber.
We will move through the fibers outside $Z$ in some order
$A_1,\dots,A_q$, where $q$ is the number of such fibers.
At each step, we will extract the prescribed powers of the variables
in the current fiber, obtaining a polynomial in the remaining
variables. Iterated extraction will give the required coefficient.

We will change the order of processing based on whether
$Z=\varnothing$ or $Z=C_p$. When $Z=\varnothing$, we will process the
fibers in the order $C_1,\dots,C_m$. When $Z=C_p$, we will process
them in the order
\[
C_1,\dots,C_{p-1},C_m,\dots,C_{p+1},
\]
leaving $Z$ unprocessed.
We will start from $P^{(0)}=P_{G,\phi}$.
Processing a fiber $A$ at step $k$ will mean choosing exponents
$d(v)$, $v\in A$, and passing from $P^{(k-1)}$ to
\[
P^{(k)}
=
\left[\prod_{v\in A}x_v^{d(v)}\right]P^{(k-1)}.
\]
An \emph{unprocessed} fiber is a fiber that has not yet been processed.
In particular, $Z$ remains unprocessed throughout.
We will orient each incomplete pair $(A,B)$ according to the
processing order, with $A$ processed before $B$, and fix all
corresponding vertex labelings independently of $\phi$.

We will see that
\begin{equation}\label{eq:invariant}
P^{(k)}
=
\pm c_k
\times \prod_{e=uv}(x_u-x_v)
\times \prod_{(A,B)}V(b_1,\dots,b_\mu),
\end{equation}
where $x_z$ is interpreted as $\phi(z)$ for $z\in Z$,
$e$ ranges over the edges with both ends unprocessed and not both
in $Z$, and $(A,B)$ ranges over the incomplete pairs with $A$
processed and $B$ unprocessed.
Here $c_k$ is the product of $\mu!$ over the incomplete pairs whose
two fibers have both been processed.
Initially, the invariant holds with $c_0=1$.

Let $\mathcal{U}_A$ be the set of unprocessed fibers other than the
current fiber $A=C_i$; this includes $Z$ when $Z$ is nonempty.
If $Z=\varnothing$, or if $Z=C_p$ and $i<p$, then
\[
\mathcal{U}_A=\{C_{i+1},\dots,C_m\}.
\]
If $Z=C_p$ and $i>p$, then
\[
\mathcal{U}_A=\{C_p,\dots,C_{i-1}\}.
\]
Thus all other unprocessed fibers lie on one side of $A$.
There are two possible scenarios.
\begin{enumerate}
\item
If the incomplete partner $B$ of $A$ is unprocessed, let $T$ be the
union of $B$ and all fibers strictly between $A$ and $B$.
These fibers are all unprocessed.
By Lemma~\ref{lem:windows}(a)--(b), applied in the appropriate
direction, $T$ is a clique containing every unprocessed neighbor
of $A$. The only missing $A$--$T$ pairs are $a_ib_i$.
Moreover, $Z\cap T=\varnothing$, since $Z$ is neither an endpoint
nor an interior fiber of an incomplete pair.

\item
Otherwise, let $T'$ be the set of unprocessed neighbors of $A$.
Since $A$ has no incomplete partner on the unprocessed side,
Lemma~\ref{lem:windows}(c), applied in the appropriate direction,
shows that $A\cup T'$ is a clique, possibly containing $Z$.
\end{enumerate}

Suppose that the invariant holds before processing $A$.
The factors of $P^{(k-1)}$ involving the variables of $A$ are exactly
the edge factors inside $A$, the edge factors from $A$ to
unprocessed fibers, and a received Vandermonde factor if the
incomplete partner of $A$ has already been processed.
By condition~\ref{condition: incomplete-partner}, at most one
received factor involves $A$, and none is present in Case~(1).
All other factors are independent of the variables of $A$.
At this step, all edge factors between $A$ and processed fibers
have disappeared from the remaining polynomial.

In Case~(1), the edge factors at $A$ form $F_A$ of \eqref{eq:FA},
and we use the exponents of Lemma~\ref{lem:sending}.
Their extraction leaves the factor
$\pm V(b_1,\dots,b_\mu)$, which will be processed when we process $B$.

In Case~(2), the edge factors at the current fiber $A$, together
with its received factor if one is present, form the receiving
polynomial $F_B$ of \eqref{eq:FB}, after relabeling.
We take $\mu=0$ when nothing is received and use the exponents of
Lemma~\ref{lem:receiving}.
The coefficient is $\pm\mu!$, which is absorbed into $c_k$.
In both cases, the invariant \eqref{eq:invariant} is maintained,
and the exponents are at most $\omega-1$ by the clique-number
bounds in the sending and receiving lemmas.
In a receiving step, vertices of $Z$ are still counted in the
clique $A\cup T'$, even though their colors occur as constants.

Since $Z$ is never an endpoint of an incomplete pair, both fibers
of every incomplete pair are processed.
When all fibers outside $Z$ have been processed, no edge factor
and no pending Vandermonde factor remains. Therefore,
\[
P^{(\mathrm{final})}
=
\pm\prod_{(U,V)}\mu(U,V)!.
\]

The exponents chosen above depend only on the graph, the elimination
order and the fixed vertex labelings, and not on $\phi$.
Specialize the coefficient identity just proved to the zero map
$\phi_0:Z\to\RR$.
The polynomial $P_{G,\phi_0}$ is homogeneous, and its degree equals
that of $P_{G,\phi}$, since both are products of the same number of
nonconstant linear factors.
The coefficient of $\prod_{v\notin Z}x_v^{d(v)}$ in
$P_{G,\phi_0}$ is still nonzero.
Homogeneity therefore gives
\[
\sum_{v\notin Z}d(v)
=
\deg P_{G,\phi_0}
=
\deg P_{G,\phi},
\]
as required.
\end{proof}

\begin{corollary}\label{cor:interval-ordinary}
Every graph $G$ of interval type is $\omega(G)$-choosable.
\end{corollary}

\begin{proof}
The empty graph is immediate.  Otherwise saturate the incomplete pairs
without changing the clique number.  Lemma~\ref{lem:elimination}, with
$Z=\varnothing$, gives a nonzero coefficient with
$d(v)+1\le\omega(G)$ at every vertex.  Corollary~\ref{cor:nullstellensatz}
colors the saturated graph from any lists of size at least $\omega(G)$.
The same coloring is proper in $G$.
\end{proof}

\subsection{Placing the prescribed clique in a central fiber}

Ordinary choosability uses full elimination.  To preserve a prescribed
coloring, we instead eliminate toward its clique.  The following two
lemmas put that clique in one fiber and ensure that no Vandermonde
transfer passes through it.  This is where simpliciality enters the
interval argument.


\begin{lemma}[Protected refinement]\label{lem:refinement}
Let $G$ be of interval type, with fixed \emph{representation}
$C_1,\dots,C_m$. Let $S$ be a clique consisting of simplicial vertices of $G$.

Then this representation can be \emph{refined} by replacing each
$C_i$ with an ordered partition of $C_i$ into one or more nonempty
cliques, placed consecutively in the position originally occupied
by $C_i$, so that the resulting ordered partition still satisfies
\ref{condition: interval} and \ref{condition: incomplete-partner}
and no new fiber meeting $S$ has an incomplete partner.
Equivalently, every new fiber meeting $S$ forms a complete or
anticomplete pair with every other new fiber.
Only the ordered partition is changed while the vertices and edges of
$G$ remain unchanged.
\end{lemma}


\begin{proof}
We will simultaneously refine all incomplete pairs having an endpoint
that meets $S$.

Consider such a pair $(U,V)$, with $U$ before $V$, and suppose first
that $T=S\cap U\ne\varnothing$.
By Lemma~\ref{lem:simplicial-clique}, the vertices of $T$ have
the same neighbors in $V$. Put
\[
U'=U\setminus T,\qquad
V_0=N(T)\cap V,\qquad
V_1=V\setminus V_0.
\]
Since a vertex of $T$ is simplicial and adjacent to every vertex
of $U'\cup V_0$, the set $V_0$ is complete to $U$.
Also, $T$ is anticomplete to $V_1$.
Moreover, $S\cap V=\varnothing$. Otherwise adjacent simplicial
vertices in $U$ and $V$ would be true twins whose common closed
neighborhood contains $U\cup V$, making $(U,V)$ complete.

We will replace $U$ by the consecutive fibers $T,U'$ and $V$ by
the consecutive fibers $V_0,V_1$, omitting empty parts.
If $S$ meets $V$ instead of $U$, we will use the reversed
construction.

These replacements are compatible.
Indeed, by \ref{condition: incomplete-partner}, distinct incomplete
pairs have no endpoint fiber in common, so each original fiber
is split at most once.
Furthermore, every complete or anticomplete relation between
old fibers is inherited by all their parts.
Within the replacement of $(U,V)$, the only possibly incomplete
pair is $(U',V_1)$, whose two fibers are disjoint from $S$.
Thus \ref{condition: incomplete-partner} is preserved, and no
new fiber meeting $S$ has an incomplete partner.

If we check \ref{condition: interval}, then we would be done.
Triples whose vertices belonged to three distinct old fibers
are covered by the original condition.
Otherwise, two vertices belong to one split old fiber.
If the third vertex belongs to an old fiber other than its
incomplete partner, that fiber is complete or anticomplete
to the split fiber, so the condition follows.

The remaining triples lie within a refined incomplete pair.
For the construction above, the parts occur in the relative order
$T,U',V_0,V_1$.
Since $T$ is anticomplete to $V_1$, the only triples of distinct
parts whose outside vertices can be adjacent have types
$(T,U',V_0)$ or $(U',V_0,V_1)$.
In both cases, the required edges exist because $V_0$ is complete
to $U$ and each old fiber is a clique.
The reversed construction satisfies the condition by the same
argument with the order reversed.

Hence the new ordered partition is a representation of interval
type, and every fiber meeting $S$ has no incomplete partner.
\end{proof}

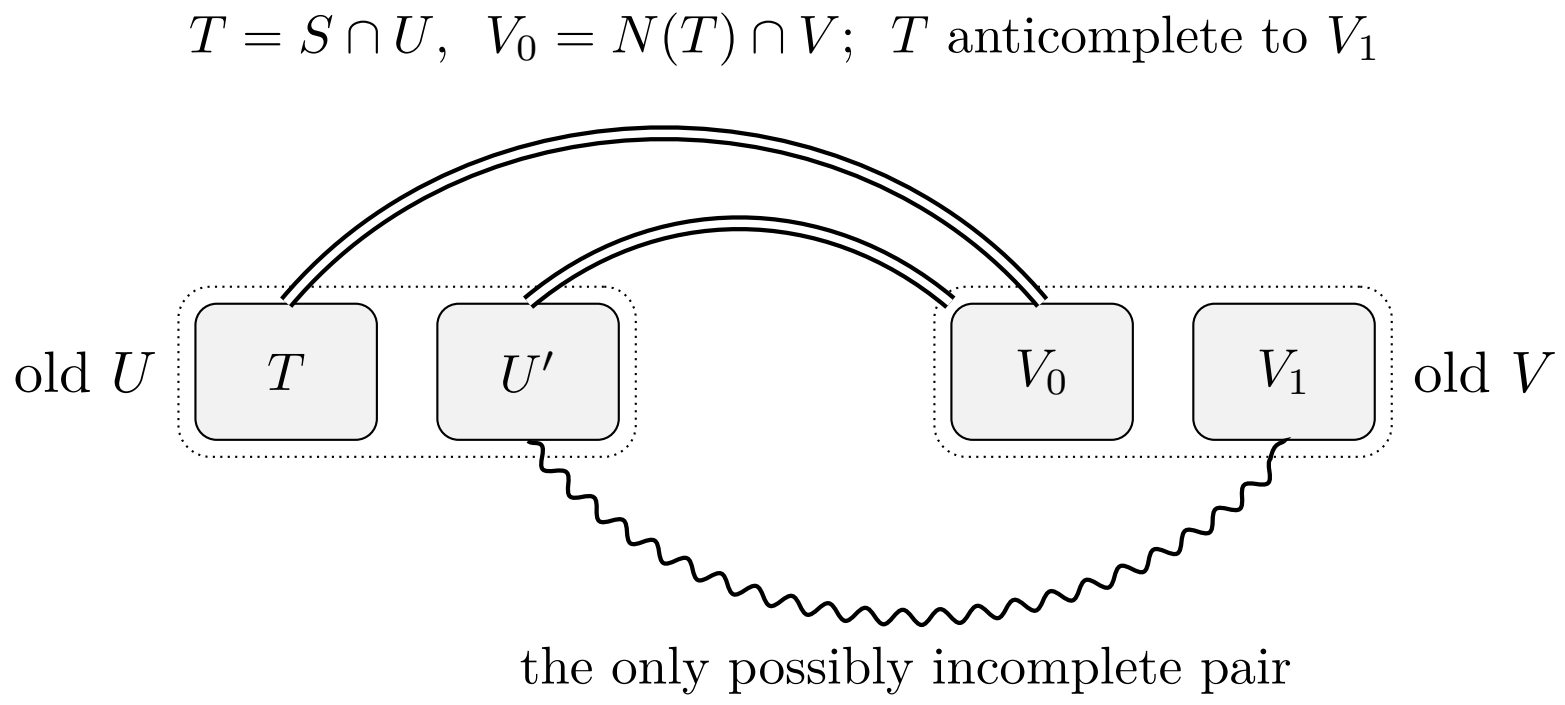
\begin{figure}[!htbp]
\centering
\begin{tikzpicture}[scale=1.0]
\node[fib,minimum width=12mm] (T) at (0,0) {$T$};
\node[fib,minimum width=12mm] (Up) at (1.6,0) {$U'$};
\node[fib,minimum width=12mm] (V0) at (5.0,0) {$V_0$};
\node[fib,minimum width=12mm] (V1) at (6.6,0) {$V_1$};
\node[draw,densely dotted,rounded corners=6pt,fit=(T)(Up),inner sep=3pt,label=left:{old $U$}] {};
\node[draw,densely dotted,rounded corners=6pt,fit=(V0)(V1),inner sep=3pt,label=right:{old $V$}] {};
\draw[full] (T.north) to[out=50,in=130] (V0.north);
\draw[full] (Up.north) to[out=40,in=140] (V0.north west);
\draw[incomp] (Up.south) to[out=-50,in=-130] node[lab,below=1mm]{the only possibly incomplete pair}(V1.south);
\node[lab] at (3.3,2.2) {$T=S\cap U$, \ $V_0=N(T)\cap V$; \ $T$ anticomplete to $V_1$};
\end{tikzpicture}
\caption{Refining an incomplete pair $(U,V)$ that meets $S$ in $T$.  Double
lines are complete pairs.  The window between $U$ and $V$ (not drawn) is
complete to all four new fibers.}
\label{fig:refinement}
\end{figure}

\begin{lemma}[Central fiber]\label{lem:central}
Let $G$ be of interval type and let $S\ne\varnothing$ be a clique
of simplicial vertices containing every simplicial vertex of $G$.  Then $G$
has an interval representation in which $S$ is a single fiber $C_p$ that has
no incomplete partner and is straddled by no incomplete pair.
\end{lemma}

\begin{proof}
    Consider the refinement produced by Lemma~\ref{lem:refinement}.
Any fiber $C$ that meets $S$ has no incomplete partner, so all
its vertices have the same closed neighborhood.
This common closed neighborhood is a clique, since $C$ contains
a simplicial vertex from $S$.
Therefore every vertex of $C$ is simplicial, and hence $C\subseteq S$.
Thus $S$ is a union of fibers.

Let $C_a$ and $C_b$ be the first and last fibers contained in $S$.
We will show that these fibers form a consecutive block.
Let $C$ be a fiber strictly between $C_a$ and $C_b$ that is not
already contained in $S$, and let $D$ be any fiber contained in $S$.
The pair $(C,D)$ cannot be incomplete, since $D$ has no incomplete
partner.
If $D$ lies before $C$, an edge between $D$ and $C_b$ spans $C$.
If $D$ lies after $C$, an edge between $C_a$ and $D$ spans $C$.
In either case, \ref{condition: interval} rules out anticompleteness.
Therefore $C$ is complete to every fiber contained in $S$.

By Lemma~\ref{lem:simplicial-clique}, the vertices of $S$
have a common closed neighborhood $K$, which is a clique.
The completeness just established gives
\[
C_a\cup\cdots\cup C_b\subseteq K.
\]
Now let $x\in C_i$ with $a<i<b$.
We will show that $N[x]\subseteq K$.
Suppose that $x$ has a neighbor $y\notin K$.
Then $y$ lies outside $C_a\cup\cdots\cup C_b$.
If $y$ lies in a fiber before $C_a$, the edge $yx$ and
\ref{condition: interval} imply that $y$ is adjacent to a vertex
of $C_a$.
If $y$ lies in a fiber after $C_b$, the same condition implies
that $y$ is adjacent to a vertex of $C_b$.
In either case, $y\in K$, a contradiction.
Thus $N[x]\subseteq K$, so $x$ is simplicial and therefore belongs
to $S$.
Consequently, every intervening fiber is contained in $S$, and
\[
S=C_a\cup\cdots\cup C_b.
\]

We will merge this consecutive block into a single fiber.
Condition~\ref{condition: interval} is preserved, since any three
vertices in distinct new fibers belonged to distinct old fibers
in the same order.
Moreover, every old fiber contained in $S$ has no incomplete
partner, and all vertices of $S$ have the same closed neighborhood.
Hence the merged fiber is complete or anticomplete to every
other fiber and has no incomplete partner.
All other fiber pairs are unchanged, so
\ref{condition: incomplete-partner} is also preserved.

Finally, suppose that an incomplete pair $(U,V)$ straddles
the merged fiber $S$.
By Lemma~\ref{lem:windows}(a), $S$ is complete to $U\cup V$.
Choose nonadjacent vertices $u\in U$ and $v\in V$.
Then every vertex of $S$ has the two nonadjacent vertices $u,v$
in its neighborhood, contradicting simpliciality.
After relabeling the fibers, $S$ is the required fiber $C_p$.
\end{proof}

\subsection{The interval precoloring theorem}

\begin{theorem}[Interval precoloring extension]\label{thm:interval-P}
Every graph $G$ of interval type satisfies \Pk\ for every $k\ge\omega(G)$.
\end{theorem}

\begin{proof}
Fix $k$ and suppose the assertion fails.  Minimize $|V(G)|$ over all
counterexamples $(G,S,\phi,L)$ in the hereditary interval-type class
with $\omega(G)\le k$.  By Lemma~\ref{lem:minimal},
$G$ is connected and every simplicial vertex of $G$ lies in $S$.  If
$S=\varnothing$, Corollary~\ref{cor:interval-ordinary} colors $G$; so
$S\ne\varnothing$, and if $V(G)=S$ there is nothing to do.

By Lemma~\ref{lem:central} we may take $S=C_p$ with no incomplete partner and
no straddling pair.  Saturate by Lemma~\ref{lem:interval-saturation}: no edge
at $S$ is added, $\omega$ is unchanged, $S$ is still a clique of simplicial
vertices, and a coloring of the supergraph extending $\phi$ is one of $G$.
Represent the finitely many colors by distinct real numbers.
Lemma~\ref{lem:elimination}, applied to the saturated graph with $Z=S$,
gives a nonzero coefficient whose exponents satisfy
$0\le d(v)\le\omega(G)-1\le k-1$ and have the required total degree.
The prescribed coloring $\phi$ is proper on $S$ by the definition of
a \Pk-instance.  Thus every hypothesis of
Corollary~\ref{cor:nullstellensatz} holds, and it extends $\phi$ from the
lists of size at least $k\ge d(v)+1$.  This contradicts the choice
of $G$.
\end{proof}
\section{Structural reductions}\label{sec:structural-reductions}

The coloring lemmas supply the required pieces of the  proof.
We first show that \Pk\ is preserved under double clique sums.
A minimum counterexample consequently has no rich clique cutset, and
deleting its prescribed simplicial clique leaves an atom.  The final
two reductions determine how this clique can attach to the two types
of atoms.

\subsection{Rich and double clique cutsets}\label{subsec:rich}

A clique cutset $K$ of a connected graph $G$ is \emph{rich} if
$V(G)\setminus K$ can be partitioned into nonempty anticomplete sets $A,B$
each of which contains two vertices that are not true twins in $G$
(the adjacent clones of~\cite{CS}).

\begin{figure}[!htbp]
\centering
\begin{tikzpicture}[scale=1.0]
\node[cbox,minimum width=16mm,minimum height=16mm] (V1) at (0,0) {$V_1$};
\node[cbox,minimum width=12mm,minimum height=12mm] (C1) at (2.4,0) {$C_1$};
\node[cbox,minimum width=12mm,minimum height=12mm] (C2) at (4.4,0) {$C_2$};
\node[cbox,minimum width=16mm,minimum height=16mm] (V2) at (6.8,0) {$V_2$};
\node[clq,fit=(C1)(C2),inner sep=3pt,label=above:{clique $D$}] {};
\draw[incomp] (V1)--(C1); \draw[full] (C1)--(C2); \draw[incomp] (C2)--(V2);
\node[lab] at (3.4,-1.6) {no edges $V_1$--$C_2$, $V_1$--$V_2$, $C_1$--$V_2$};
\end{tikzpicture}
\caption{A double clique cutset.  Wavy lines denote arbitrary adjacency.}
\label{fig:double}
\end{figure}
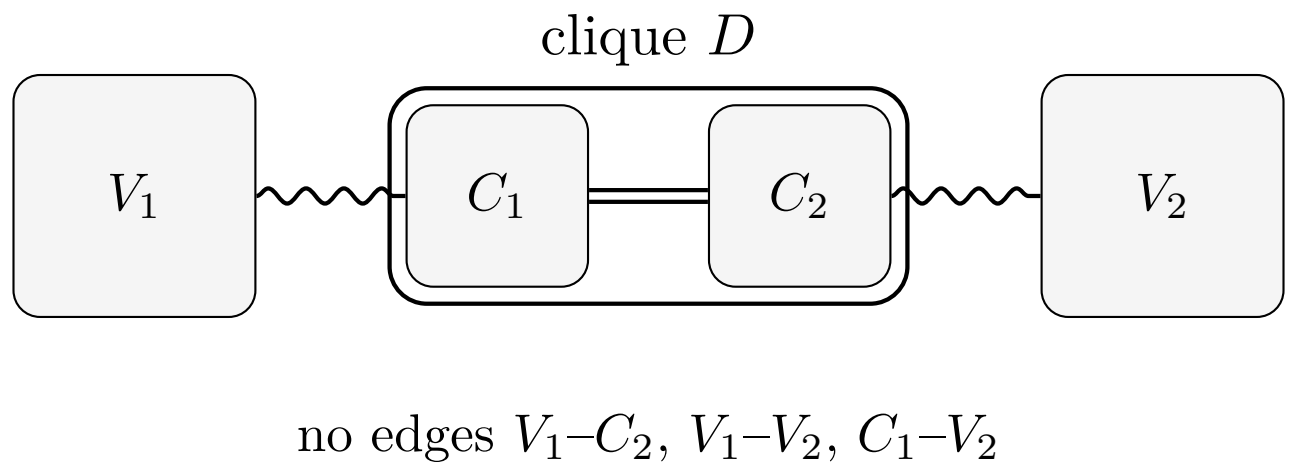

A \emph{double clique cutset} of a graph $G$ is a clique cutset
$D=C_1\dcup C_2$ for which
\[
V(G)=V_1\dcup C_1\dcup C_2\dcup V_2,
\qquad V_1,V_2\ne\varnothing,
\]
with $V_1$ anticomplete to $C_2\cup V_2$ and $V_2$ anticomplete to $C_1$.
In this case $G$ is the \emph{double clique-sum} of $G[V_1\cup D]$ and
$G[V_2\cup D]$ (Figure~\ref{fig:double}).

The next lemma shows that \Pk\ is preserved by double clique sums.

\begin{lemma}[Double clique sums]\label{lem:double-P}
Let $G$ be the double clique-sum of $G_1=G[V_1\cup D]$ and
$G_2=G[V_2\cup D]$, where $D=C_1\dcup C_2$, and let $k\ge\omega(G)$.
If both $G_1$ and $G_2$ satisfy \Pk, then $G$ satisfies \Pk.
\end{lemma}

\begin{proof}
Let $(G,S,\phi,L)$ be a \Pk-instance. Since $S$ is a clique and $V_1$
is anticomplete to $V_2$, we may assume by symmetry that
$S\cap V_2=\varnothing$. Thus $S\subseteq V(G_1)$, and its vertices
remain simplicial in $G_1$. Apply \Pk\ in $G_1$ to obtain a coloring
$f_1$ extending $\phi$ and respecting the lists outside $S$.

In $G_2$, the set
\[
T=C_1\cup(S\cap C_2)
\]
is a clique of simplicial vertices: every vertex of $C_1$ has all its
neighbors in the clique $D$, and vertices of $S\cap C_2$ remain
simplicial. Precolor $T$ by $f_1|_T$. This coloring agrees with $\phi$
on $S\cap C_2$, and every vertex of $G_2-T$ has its original list of
at least $k$ colors. Apply \Pk\ in $G_2$ to obtain a coloring $f_2$
extending $f_1|_T$ and respecting these lists.

Define $f=f_1$ on $V_1\cup C_1$ and $f=f_2$ on $V_2\cup C_2$.
The only edges between these two parts are the $C_1$--$C_2$ edges;
$f_2$ is proper on them and agrees with $f_1$ on $C_1$. Thus $f$ is
proper, respects the lists outside $S$, and extends $\phi$.
\end{proof}

We recall the following structural result of Chudnovsky and Seymour.

\begin{lemma}[Chudnovsky--Seymour {\cite[Theorem~1.4]{CS}}]\label{thm:rich-structure}
A connected claw-free graph with a rich clique cutset is of interval type or
admits a double clique cutset.
\end{lemma}

Consequently, by Lemmas~\ref{thm:rich-structure} and~\ref{lem:double-P}
and Theorem~\ref{thm:interval-P}, a minimum-order counterexample to
\Pk\ in the class of claw-free perfect graphs has no rich clique cutset.

\subsection{Atoms after removing the prescribed clique}

A noncomplete atom has no simplicial vertex: a nonuniversal simplicial vertex
$v$ has the clique cutset $N(v)$, and a universal simplicial vertex makes the
graph complete.

\begin{lemma}\label{lem:no-rich-atom}
Let $G$ be a connected graph with no rich clique cutset, and let $S$ be
a clique of simplicial vertices containing every simplicial vertex of
$G$. Then $G-S$ is either empty or an atom.
\end{lemma}

\begin{proof}
Put $H=G-S$ and $Q=N_G(S)$, and suppose $H\ne\varnothing$.
If $S=\varnothing$, then $H=G$ is connected. Otherwise, every component
of $H$ contains a vertex of $Q$, since $G$ is connected. As $Q$ is a
clique by Lemma~\ref{lem:simplicial-clique}, $H$ is connected in this
case as well.

Suppose that $D$ is a clique cutset of $H$, and partition
$V(H)\setminus D$ into nonempty anticomplete sets $V_1,V_2$.
Since $Q$ is a clique, we may assume that $Q\setminus D\subseteq V_1$.
Then $V_1\cup S$ and $V_2$ are nonempty anticomplete sets partitioning
$V(G)\setminus D$.

We show that each side contains two vertices that are not true twins
in $G$. If $V_2$ contained no such pair, it would be a clique of
pairwise true twins (or a singleton). Each vertex of $V_2$ has
neighbors only in $V_2\cup D$; its neighbors in $D$ are complete to
$V_2$, by the twin property. Since $D$ is a clique, every vertex of
$V_2$ would be simplicial in $G$, contrary to $V_2\cap S=\varnothing$.
If $S=\varnothing$, the same argument applies to $V_1$.
If $S\ne\varnothing$, choose $s\in S$ and $v\in V_1$. They cannot be
true twins, since that would make $v$ simplicial as well, contrary to
$v\notin S$.

Thus $D$ is a rich clique cutset of $G$, a contradiction. Hence $H$
has no clique cutset and is an atom.
\end{proof}

\subsection{Nonempty attachments to peculiar graphs}\label{sec:peculiar}

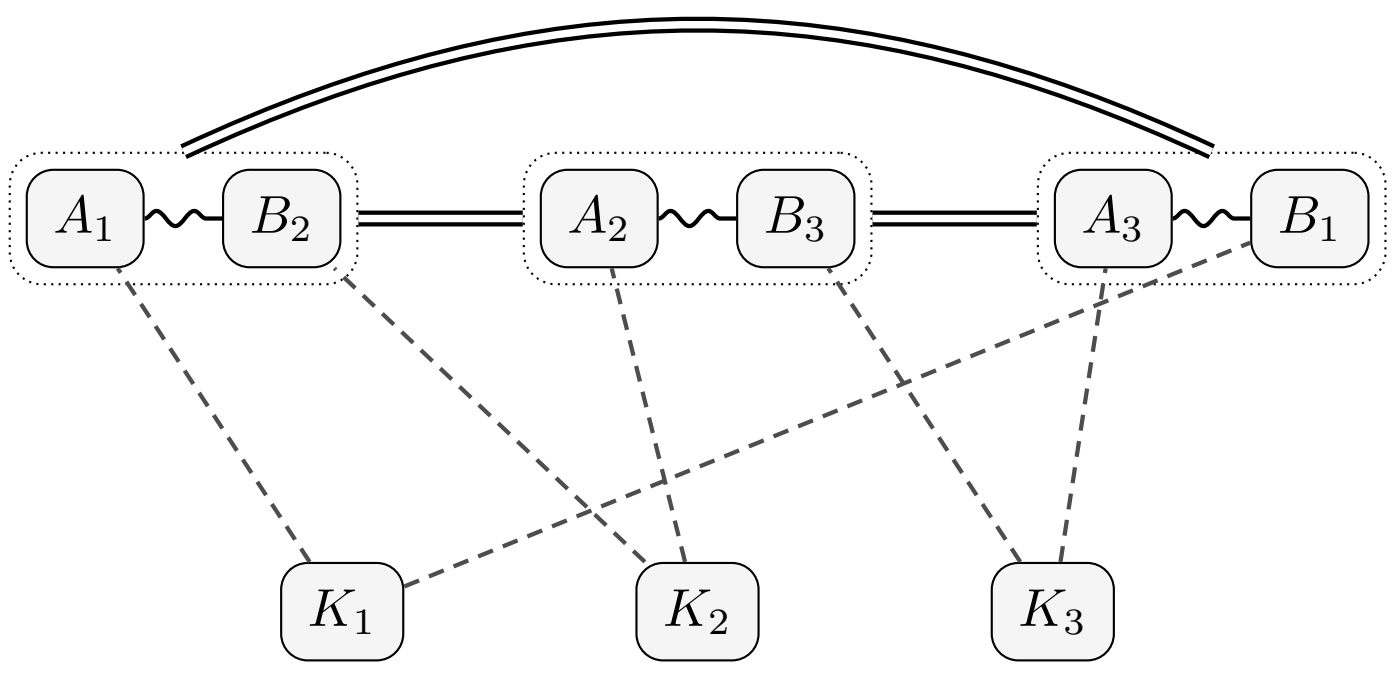
\begin{figure}[!htbp]
\centering
\begin{tikzpicture}[scale=1.0]
\node[cbox] (A1) at (0,0) {$A_1$}; \node[cbox] (B2) at (1.3,0) {$B_2$};
\node[cbox] (A2) at (3.4,0) {$A_2$}; \node[cbox] (B3) at (4.7,0) {$B_3$};
\node[cbox] (A3) at (6.8,0) {$A_3$}; \node[cbox] (B1) at (8.1,0) {$B_1$};
\draw[incomp] (A1)--(B2); \draw[incomp] (A2)--(B3); \draw[incomp] (A3)--(B1);
\node[draw,densely dotted,rounded corners=6pt,fit=(A1)(B2),inner sep=3pt] (P1) {};
\node[draw,densely dotted,rounded corners=6pt,fit=(A2)(B3),inner sep=3pt] (P2) {};
\node[draw,densely dotted,rounded corners=6pt,fit=(A3)(B1),inner sep=3pt] (P3) {};
\draw[full] (P1)--(P2); \draw[full] (P2)--(P3);
\draw[full] (P1.north) to[bend left=25] (P3.north);
\node[cbox] (K1) at (1.7,-2.6) {$K_1$};
\node[cbox] (K2) at (4.05,-2.6) {$K_2$};
\node[cbox] (K3) at (6.4,-2.6) {$K_3$};
\draw[nonedge] (K1)--(A1); \draw[nonedge] (K1)--(B1);
\draw[nonedge] (K2)--(A2); \draw[nonedge] (K2)--(B2);
\draw[nonedge] (K3)--(A3); \draw[nonedge] (K3)--(B3);
\node[lab,align=center,text width=12cm] at (4.05,-3.85)
{Dashed lines show the anticomplete pairs $K_i$--$A_i$ and $K_i$--$B_i$.\\
Every other $K_i$--core pair is complete.};
\end{tikzpicture}
\caption{A peculiar graph. Each dotted box is a cobipartite piece with at
least one nonadjacent pair; its wavy line represents a noncomplete
adjacency, possibly anticomplete. Distinct pieces are complete to each
other, and the $K_i$ are pairwise anticomplete.}
\label{fig:peculiar}
\end{figure}

\begin{definition}[Peculiar graphs {\cite{ChvatalSbihi}}]
A graph $G$ is called \emph{peculiar} if its vertex set admits a
partition into nine pairwise disjoint cliques
$(A_1,B_1,A_2,B_2,A_3,B_3,K_1,K_2,K_3)$ satisfying the following
conditions, with indices taken modulo $3$:
\begin{enumerate}
\item The six cliques $A_i,B_i$ are nonempty. We allow the $K_i$ to
be empty; the standard nondegenerate definition requires all nine
cliques to be nonempty.
\item Each pair $(A_i,B_{i+1})$ contains a nonadjacent pair, and the
three sets $A_1\cup B_2$, $A_2\cup B_3$, $A_3\cup B_1$ are pairwise
complete to one another.
\item The cliques $K_1,K_2,K_3$ are pairwise anticomplete.
\item For each $i$, the clique $K_i$ is anticomplete to $A_i\cup B_i$
and complete to $A_j\cup B_j$ for every $j\ne i$.
\end{enumerate}
See Figure~\ref{fig:peculiar}. We call the $A_i,B_i$ the
\emph{core cliques} and the $K_i$ the \emph{non-core cliques}.
\end{definition}

\begin{lemma}\label{lem:peculiar-no-attachment}
A peculiar graph is either cobipartite or has no nonempty attachable
clique.
\end{lemma}

\begin{proof}
If a non-core clique is empty, say $K_3=\varnothing$, then
\[
K_1\cup A_2\cup B_2\cup A_3,
\qquad K_2\cup A_1\cup B_1\cup B_3
\]
are cliques partitioning the vertex set. Thus the graph is cobipartite.
Assume that all three $K_i$ are nonempty, and suppose that $Q$ is a
nonempty attachable clique. 

 By  \eqref{eq:att-rule}, two distinct nonadjacent neighbors of $q$ cannot both lie
outside $Q$.

Suppose first that $Q$ meets a core clique. The cyclic shift of indices
and the map $A_i\mapsto B_{-i}$, $B_i\mapsto A_{-i}$,
$K_i\mapsto K_{-i}$ preserve the definition and act transitively on
the core cliques. Thus we may assume that $a\in A_1\cap Q$.
Since $a$ is complete to the nonempty anticomplete cliques $K_2,K_3$,
attachability implies that one of them lies entirely in $Q$; the other
is disjoint from $Q$, since $Q$ is a clique. Also $K_1\cap Q=\varnothing$.

If $K_2\subseteq Q$, a vertex of $K_3$ is an outside neighbor of $a$
nonadjacent to every vertex of $A_3\cup B_3$. These vertices are all
adjacent to $a$, so \eqref{eq:att-rule} gives $A_3\cup B_3\subseteq Q$.
For $b_3\in B_3$, an outside neighbor in $K_1$ similarly forces
$B_1\subseteq Q$. Then $Q$ contains $A_3\cup B_1$, which has a
nonadjacent pair, a contradiction.

If $K_3\subseteq Q$, outside neighbors in $K_2$ first force
$A_2\subseteq Q$. For $a_2\in A_2$, outside neighbors in $K_1$ force
$A_1\cup B_1\subseteq Q$. For $b_1\in B_1$, outside neighbors in
$K_2$ then force $B_2\subseteq Q$. Now $Q$ contains the noncomplete
pair $(A_1,B_2)$, again a contradiction.

Finally, if $Q$ meets no core clique, then $Q\subseteq K_i$ for some
$i$. A vertex of $Q$ is complete to $A_{i+1}\cup B_{i+2}$, which
contains a nonadjacent pair outside $Q$. This contradicts attachability.
\end{proof}

\subsection{Attaching a simplicial vertex to an elementary atom}
\label{sec:elementary-atom-attachment}

In this subsection we show that adjoining a prescribed simplicial clique
to a non-cobipartite elementary atom gives an elementary graph
(Corollary~\ref{cor:elementary-atom-twin-extension}). By Lemma~\ref{lem:elementary-twins} it suffices to add a single vertex $s$
whose neighbourhood is an attachable clique $Q$ of the atom $H$. The
remaining vertices of the clique are then added as true twins of $s$.
We work directly with the defining edge $2$-colouring rather than with
the root. First we show that $Q$ is either an entire \emph{hub clique}
$C_u$ (the clique formed by the vertices coming from the root edges at a
root vertex $u$), or meets both shores of a single augment in a set
$X_0 \cup Y_0$ of a restricted form. In either case, every edge from $Q$
to the rest of $H$ receives the same colour. Giving each new edge $sq$
the other colour then makes every new induced path $s$--$q$--$x$
bichromatic.

\begin{lemma}\label{lem:elementary-twins}
Adding a true twin of a vertex to an elementary graph gives an elementary
graph.
\end{lemma}

\begin{proof}
Start with an elementary two-coloring of the edges of the old graph.
When $v$ is added as a true twin of $u$, color $vx$ as $ux$ for every
old neighbor $x$ of $u$, and color $uv$ arbitrarily.  The edge $uv$ lies in no induced three-vertex path, and
every other induced path $x\,v\,y$ or $x\,y\,v$ corresponds to
$x\,u\,y$ or $x\,y\,u$ in the old graph.
\end{proof}

\begin{lemma}\label{lem:elementary-atom-extension}
Let $H$ be an elementary atom that is not cobipartite, and let $Q$ be an
attachable clique of $H$.  Adding a vertex $s$ with $N(s)=Q$ gives an
elementary graph.
\end{lemma}

\begin{proof}[Proof:]
If $Q=\varnothing$, the new vertex is isolated and the assertion is
immediate. Assume $Q\ne\varnothing$. By the Maffray--Reed structure
theorem, fix a normalized representation of $H$ with bipartite root $B$,
and absorb every complete augment by Lemma~\ref{lem:augment}(c).
Thus every remaining augment is noncomplete. Since $H$ is a
noncomplete atom, it has no simplicial vertex.

For each augment $(X,Y)$, both exterior neighborhoods $N_X,N_Y$ are
nonempty: otherwise deleting the opposite shore would separate a shore
from a vertex outside $X\cup Y$, which exists because $H$ is not
cobipartite. Every vertex of either shore has a neighbor in the other,
since otherwise it would be simplicial.

Every nonisolated root vertex has at least two distinct neighbors.
Indeed, if $u$ had only neighbor $v$, a root edge $e=uv$ would be
simplicial in $\mathcal L(B)$. If $e$ were unaugmented, its neighbors
would still form a clique after augmentation: no flat pair can lie
entirely among them, since $e$ would complete a triangle. If $e$
belonged to a designated flat pair $ef$, simpliciality and flatness
would make $f$ its only neighbor, giving an augment shore with empty
exterior neighborhood. Both possibilities are excluded.

We use the elementary $2$-edge-coloring arising from this representation.

Fix a bipartition $(A,B')$ of the root multigraph. If two nonparallel
root edges meet at $u$, color their line-graph edge according to the
part containing $u$. Edges joining vertices represented by parallel
root edges may be colored arbitrarily: those vertices are true twins
and their common edge lies in no induced three-vertex path. For an
augment replacing a flat edge corresponding to two
root edges $uc,\;cv$
where $d(c)=2$, let its two shores be $X$ and $Y$. Since $u$ and $v$
belong to the same side of the bipartition, say with colour $\alpha$,
colour all edges inside $X$ and inside $Y$, and all edges joining either
shore to its exterior attachment clique, with colour $\alpha$; colour
the edges between $X$ and $Y$ with the other colour. This coloring is
elementary: an induced three-vertex path that uses an edge within a
shore or across an augment has one shore-side edge and one cross-shore
edge; every other such path projects to an induced three-vertex path
in the original line graph, whose two edges receive different colors.

For a nonisolated root vertex $u$ that is not the center of a
designated augmented pair, let $C_u$ denote the corresponding
\emph{hub clique}: the vertices arising from root edges incident with
$u$, with an augmented endpoint replaced by its shore. Its incidence
blocks are indexed by the distinct neighbors of $u$: a parallel bundle
forms one block, and an augmented root edge contributes its entire
shore. The set $C_u$ is a clique, since no designated flat pair has
both root edges incident with $u$.

We first record the form of an attachable clique $Q$.

\medskip
\noindent
\emph{Claim.}
Either
\begin{enumerate}
    \item $Q=C_u$ for some hub clique $C_u$, or
    \item $Q$ meets both shores $X,Y$ of a single augment and
    \[
        Q=X_0\cup Y_0,
        \qquad
        \varnothing\neq X_0\subseteq X,\quad
        \varnothing\neq Y_0\subseteq Y,
    \]
    where $X_0$ is complete to $Y_0$, $X_0$ is anticomplete to
    $Y\setminus Y_0$, and $Y_0$ is anticomplete to $X\setminus X_0$.
\end{enumerate}

Suppose first that $Q$ is contained in a hub clique $C_u$. The root
degree observation gives at least two incidence blocks. Every vertex
of $C_u$ has a neighbor outside $C_u$, since otherwise it would be
simplicial. If $q$ belongs to the block indexed by a neighbor $v$ of
$u$, an outside neighbor of $q$ arises through $v$, or in the opposite
shore of its augment. Such a neighbor is nonadjacent to every other
block of $C_u$, by bipartiteness and Lemma~\ref{lem:augment}(a).

If $Q\neq C_u$, choose $q\in Q$ and $r\in C_u\setminus Q$ from different
incidence blocks; this is possible because there are at least two
blocks and both $Q$ and $C_u\setminus Q$ are nonempty. Choose a
neighbor $z\notin C_u$ of $q$ through the block containing $q$. By the augmentation geometry, $r$ is not adjacent
to $z$. Hence
\[
        r,z\in N_H(q)\setminus Q
\]
are nonadjacent, contradicting the assumption that $Q$ is attachable.
Thus $Q=C_u$.

Now suppose that $Q$ is not contained in a hub clique. The root
edges whose corresponding vertices or shores meet $Q$ are pairwise
incident, hence share a common endpoint, since $B$ is bipartite. This
endpoint must be the center of a designated augmented pair; otherwise
$Q$ would lie in its hub. Therefore $Q$ meets both shores $X,Y$ of
that augment (a single shore lies in its exterior hub). No outside
vertex sees both shores, by Lemma~\ref{lem:augment}(a), so
\[
        Q\subseteq X\cup Y.
\]
Put
\[
        X_0=Q\cap X,\qquad Y_0=Q\cap Y.
\]
Both sets are nonempty, and, since $Q$ is a clique, $X_0$ is complete to
$Y_0$.

Let $C_u$ and $C_v$ be the two exterior hub cliques of this augment.
Since $H$ is a non-cobipartite atom, each of these hub cliques has another
incidence block. Choose
\[
        a\in C_u\setminus X,\qquad b\in C_v\setminus Y.
\]
If some $x\in X_0$ had a neighbour $y\in Y\setminus Y_0$, then
\[
        a,y\in N_H(x)\setminus Q,
\]
whereas $a$ and $y$ are nonadjacent. This contradicts attachability.
Hence $X_0$ is anticomplete to $Y\setminus Y_0$. By symmetry,
$Y_0$ is anticomplete to $X\setminus X_0$. This proves the claim.

We now add the new vertex $s$ with $N(s)=Q$ and extend the elementary
edge-colouring.

Suppose first that $Q=C_u$. Let $\alpha$ be the colour corresponding to
the root vertex $u$. Every edge
\[
        qx,\qquad q\in Q,\quad x\in N_H(q)\setminus Q,
\]
has the colour opposite to $\alpha$: such an edge leaves the hub $C_u$
through the other endpoint of the corresponding root edge, or through
the opposite shore of an augment. Colour every edge $sq$, $q\in Q$,
with colour $\alpha$.

Suppose next that
\[
        Q=X_0\cup Y_0
\]
is as in the second case of the claim. Let the corresponding flat edge
come from the root edges $uc$ and $cv$, and let $\alpha$ be the common
colour associated with $u$ and $v$. By the claim, a vertex of $Q$ has no
cross-shore neighbour outside $Q$. Hence every edge joining a vertex of
$Q$ to a vertex of $H-Q$ is a hub-side edge and therefore has colour
$\alpha$. Colour every edge $sq$, $q\in Q$, with the other colour.

In either case, the only new induced paths on three vertices are of the
form
\[
        s-q-x,
        \qquad
        q\in Q,\quad x\in N_H(q)\setminus Q.
\]
Indeed, no path $q-s-q'$ with $q,q'\in Q$ is induced, since $Q$ is a
clique. By construction, the two edges of every new induced $P_3$ have
different colours. All induced $P_3$'s contained in $H$ were already
properly bichromatic.

Thus the extended graph admits an elementary edge-colouring, and hence
the graph obtained from $H$ by adding $s$ with $N(s)=Q$ is elementary.
\end{proof}

\begin{corollary}\label{cor:elementary-atom-twin-extension}
Let $G$ be connected and claw-free, let $S\ne\varnothing$ be a clique of
simplicial vertices, and suppose $G-S$ is an elementary atom that is not
cobipartite.  Then $G$ is elementary.
\end{corollary}

\begin{proof}
By Lemma~\ref{lem:simplicial-clique}, $Q=N(S)$ is an attachable clique of
$G-S$ and every vertex of $S$ has neighborhood $Q\cup(S\setminus\{s\})$.  Add
one vertex of $S$ by Lemma~\ref{lem:elementary-atom-extension} and the others
as true twins by Lemma~\ref{lem:elementary-twins}.
\end{proof}

\section{Proof of the main theorem}\label{sec:global-P}

Gravier and Maffray~\cite[Corollary~2, p.~89]{GravierMaffray}
proved that every peculiar graph $H$ satisfies
$\ch(H)=\chi(H)$. Their cliques $Q_i$ correspond to the
cliques denoted here by $K_i$. In our application, $H$
is an induced subgraph of the perfect graph $G$, so
\[
  \ch(H)=\chi(H)=\omega(H)\le\omega(G)\le k.
\]

Hence $H$ is colorable from any lists of size at least $k$.
In the proof below, we invoke this result only in the branch
$S=\varnothing$, where no colors are prescribed, and only when
$K_1,K_2,K_3$ are all nonempty. If some $K_i$ is empty,
then $H$ is cobipartite and Theorem~\ref{thm:M} applies instead.

\begin{proof}[Proof of Theorem~\ref{thm:main-P}]
Fix $k$ and suppose the theorem fails.  Claw-free perfect graphs form an
induced-subgraph-closed class, so Lemma~\ref{lem:minimal} provides a
counterexample $(G,S,\phi,L)$ of minimum order in which $G$ is connected,
every simplicial vertex lies in $S$, and every proper induced subgraph
satisfies \Pk.

$G$ has no double clique cutset, since its pieces would be proper induced
subgraphs and Lemma~\ref{lem:double-P} would color $G$; and $G$ is not of
interval type by Theorem~\ref{thm:interval-P}.  Hence, by
Lemma~\ref{thm:rich-structure}, $G$ has no rich clique cutset.

If $V(G)=S$ there is nothing to prove.  Otherwise $H=G-S$ is a nonempty atom
by Lemma~\ref{lem:no-rich-atom}, hence elementary or peculiar by
Theorem~\ref{thm:atom}.

\emph{Empty prescribed clique.}  If $S=\varnothing$, then $G=H$ is colored by Theorem~\ref{thm:elementary-P}
if elementary, by Theorem~\ref{thm:M} (with $C=\varnothing$) if cobipartite,
and by Gravier--Maffray \cite{GravierMaffray} if peculiar with
all $K_i$ nonempty.  So $S\ne\varnothing$.

\emph{Nonempty prescribed clique.}  Since $H\ne\varnothing$ and $G$
is connected, $Q=N(S)$ is nonempty.  Lemma~\ref{lem:simplicial-clique}
makes it an attachable clique of $H$, and deleting $\phi(S)$ from the lists on
$Q$ gives an \Mk-instance $(H,Q,|S|,L')$ whose colorings extend $\phi$.
\begin{itemize}
\item If $H$ is cobipartite, Corollary~\ref{cor:cobipartite-Mk} colors this
instance.
\item If $H$ is peculiar and not cobipartite, all $K_i$ are nonempty.
Lemma~\ref{lem:peculiar-no-attachment} excludes the \emph{nonempty}
attachable clique $Q$.
\item If $H$ is elementary and not cobipartite,
Corollary~\ref{cor:elementary-atom-twin-extension} shows that $G$ itself is
elementary, and Theorem~\ref{thm:elementary-P} extends $\phi$.
\end{itemize}
In every case $\phi$ extends, contradicting the choice of $G$.
\end{proof}

\begin{proof}[Proof of Theorem~\ref{thm:main}]
Theorem~\ref{thm:main-P} with $S=\varnothing$ and $k=\omega(G)$ gives
$\ch(G)\le\omega(G)$, and $\omega(G)\le\chi(G)\le\ch(G)$ always.
\end{proof}

\section{Consequences}\label{sec:consequences}

\begin{corollary}\label{cor:global-M}
Every claw-free perfect graph $G$ satisfies \Mk\ for every $k\ge\omega(G)$.
\end{corollary}

\begin{proof}
Let $(G,Q,r,L)$ be an \Mk-instance; if $r=0$ apply Theorem~\ref{thm:main-P}
with $S=\varnothing$.  For $r>0$, form $G^+$ by adding a clique $S$ of $r$ new vertices,
complete to $Q$ and anticomplete to $V(G)\setminus Q$.

We check the hypotheses needed for \Pk.  Each new vertex is simplicial.
For $q\in Q$, its neighborhood in $G^+$ is covered by the two cliques
$(Q\setminus\{q\})\cup S$ and $N_G(q)\setminus Q$.
At every other old vertex the neighborhood is unchanged.
Hence $G^+$ is claw-free.

To check perfection, consider any induced subgraph $J$ of $G^+$.
The graph $J-S$ is perfect and can be colored with at most $\omega(J)$
colors.  Insert the vertices of $J\cap S$ one at a time.  Each is
simplicial in $J$ and has at most $\omega(J)-1$ neighbors, so an available
color exists.  Thus $\chi(J)\le\omega(J)$, proving that $G^+$ is perfect.
Finally, every clique meeting $S$ lies in $Q\cup S$, and the mixed-list
hypothesis gives $|Q|+r\le k$.  Consequently $\omega(G^+)\le k$.

Choose $r$ distinct colors absent from all original lists, precolor
$S$ with them, and add them to every list on $Q$.  The enlarged lists
on $Q$ have size at least $(k-r)+r=k$, and all other old lists already
have at least $k$ colors.  Theorem~\ref{thm:main-P} extends this
precoloring; $Q$ cannot use the fresh colors, being complete to $S$, so the
restriction to $G$ is an $L$-coloring.
\end{proof}

The following corollary recovers the theorem of Peterson and Woodall
\cite{PetersonWoodall,PetersonWoodallErratum}.

\begin{corollary}\label{cor:perfect-line}
Every perfect line graph is chromatic-choosable.  Equivalently, if the line
graph of a loopless multigraph $B$ is perfect, then $\ch'(B)=\chi'(B)$.
\end{corollary}

\begin{proof}
Line graphs are claw-free: the neighbors of the vertex $uv$ of $\mathcal L(B)$ are
covered by the two cliques of edges at $u$ and at $v$.  Apply
Theorem~\ref{thm:main}; proper vertex colorings of $\mathcal L(B)$ are proper edge
colorings of $B$, with the same correspondence for lists.
\end{proof}




\subsection*{Acknowledgements}
Claude and GPT-6 were used to assist in developing and formalizing some
of the arguments. The proofs have been written by hand and verified by the authors. 


\appendix

\section{Notation table}

\begingroup
\small
\setlength{\LTpre}{.6\baselineskip}
\setlength{\LTpost}{.6\baselineskip}
\renewcommand{\arraystretch}{1.12}
\begin{longtable}{@{}p{.26\textwidth}p{.69\textwidth}@{}}
\caption{Notation and scope. The references in the right column give the
formal definitions of specialized objects.}\label{tab:notation}\label{tab:phase2-notation}\\
\toprule
Notation & Meaning and scope\\
\midrule
\endfirsthead
\multicolumn{2}{@{}l}{\textit{Table~\thetable\ (continued)}}\\
\toprule
Notation & Meaning and scope\\
\midrule
\endhead
\midrule
\multicolumn{2}{r@{}}{\textit{Continued on the next page}}\\
\endfoot
\bottomrule
\endlastfoot
\multicolumn{2}{@{}l}{\textit{Graphs, lists, and extension properties}}\\
$G[A]$, $G-A$ & Induced subgraphs on $A$ and $V(G)\setminus A$.\\
$N_G(v)$, $N_G[v]$ & Open and closed neighborhoods.\\
$N_G(A)$ & External neighborhood of $A$, excluding $A$ itself.\\
$\omega(G),\chi(G),\ch(G)$ & Clique number, chromatic number, and choice number.\\
$L(v)$, $L(T)$ & The list at $v$ and $\bigcup_{v\in T}L(v)$.\\
$f(T)$, $\phi$ & Colors used on $T$; a prescribed proper coloring.\\
$S$, $Q$, $r$, $k$ & A clique of simplicial vertices; an attachable clique;
a nonnegative list deficit; a fixed integer list bound.\\
\Pk, \Mk & Precoloring extension and mixed-list properties
(Definition~\ref{def:PM}); in \Mk, $|Q|+r\le k$.\\
$\mathcal C$ & A class closed under induced subgraphs; fixed before
choosing a minimum counterexample.\\
\addlinespace
\multicolumn{2}{@{}l}{\textit{Nonedges and polynomial coefficients}}\\
$F$, $\nu(F)$ & The bipartite nonedge graph of a cobipartite graph and
its maximum matching size (Sections~\ref{sec:prelim}--\ref{sec:cobipartite}).\\
$K_W-M$ & The complete graph on $W$ with matching pairs $M$ deleted.\\
$P_G$, $[\mathbf{x}^d]P$ & Graph polynomial and coefficient of the indicated monomial.\\
$Z$, $P_{G,\phi}$ & Fixed-color set and polynomial with its internal edges
omitted (Corollary~\ref{cor:nullstellensatz}); $\phi$ is proper when the
polynomial is used for coloring extension.\\
$\Van(\zeta_1,\ldots,\zeta_h)$ & Vandermonde product
$\prod_{i<j}(\zeta_j-\zeta_i)$, equal to $1$ for $h\le1$.\\
$m,t,u$; $q_0,q$ & Canonical-graph parameters:
$m,t\ge1$, $u\ge0$, $q_0=m+t$, $q=m+t+u$
(Section~\ref{sec:cobipartite}).\\
\addlinespace
\multicolumn{2}{@{}l}{\textit{Line graphs and elementary graphs}}\\
$B$, $\mathcal L(B)$ & A loopless bipartite root multigraph and its line graph.\\
$\st_B(w)$, $E(w,u)$ & Root edges at $w$; the parallel bundle joining $w,u$.\\
$\Delta(B)$, $d_D^+(v)$ & Maximum root degree, counting multiplicities;
number of distinct out-neighbors in a orientation $D$.\\
$(X,Y)$, $N_X,N_Y$ & An augment and its external neighborhoods
(Lemma~\ref{lem:augment}).\\
$A_i,Q_i,M_i,A_i^*$ & Unprecolored side, its external neighborhood,
matched vertices, and unmatched vertices (Lemma~\ref{lem:GP9}).\\
$a_i,\alpha_i,\mu$ & Sizes $|A_i|,|A_i^*|,|M|$ in that lemma.\\
$m_i$, $R(v)$ & Minimum list-union excess on $A_i^*$;
residual list at $v$ in that lemma.\\
\addlinespace
\multicolumn{2}{@{}l}{\textit{Interval representations and the structural reduction}}\\
$C_1,\ldots,C_m$ & Ordered clique fibers of an interval representation
(Definition~\ref{def:interval}).\\
$W$, $\mu$ & Window of an incomplete pair; size of its matching of nonedges
(Section~\ref{sec:interval}).\\
$R$, $X_e$, $F$ & Root multigraph, clique represented by root edge $e$,
and selected root edges (Lemma~\ref{lem:elementary-atom-extension}).\\
$A_i,B_i,K_i$ & The nine cliques of a peculiar representation,
with indices modulo $3$ (Section~\ref{sec:peculiar}).\\
\end{longtable}
\endgroup
\clearpage

\end{document}